\documentclass{Yama-axSIH}

\newtheorem{theorem}{Theorem}

\newtheorem{lemma}[theorem]{Lemma}

\theoremstyle{definition}

\newtheorem{definition}[theorem]{\font\=cmssi10\Definition\bf}
\newtheorem{definitions}[theorem]{\font\=cmssi10\Definitions\bf}
\newtheorem{example}[theorem]{\font\=cmssi10\Example\bf}

\newtheorem{remark}[theorem]{\font\=cmssi10\Remark\bf}
\newtheorem{remarks}[theorem]{\font\=cmssi10\Remarks\bf}

\newcommand\Cal{\mathcal}
\newcommand\bosy{\boldsymbol}

\def\Alf{\hbox{\font\=cmmi10 scaled\magstep1\\char'013}\kern0.15mm}
\def\Eps{\hbox{\font\=cmmi10 scaled\magstep1\\char'017}\kern0.15mm}
\def\Iota{\kern.15mm\hbox{\font\=cmmi10 scaled\magstep1\\char'023}\kern0.2mm}
\def\Nu{\hbox{\font\=cmmi10 scaled\magstep1\\char'027}\kern0.25mm}
\def\uvarPi{\kern.15mm\underline{\kern-.15mm\varPi\kern-.85mm}\kern.85mm}
\def\uOmega{\kern.3mm\underline{\kern-.3mm\Omega\kern-.3mm}\kern.3mm}

\def\tauu{\char'034}
\def\fRe{\hbox{\font\=cmr9\f\kern.1mm}\roman{Re}\kern.75mm}
\def\fIm{\hbox{\font\=cmr9\f\kern.1mm}\roman{Im}\kern.65mm}
\def\vecc#1{\kern-.5mm\vec{\kern.5mm#1}}
\def\TVS{\roman{TVS}\kern0.37mm}%
\def\LCS{\roman{LCS}\kern0.37mm}%
\def\BaS{\roman{BaS}\kern0.37mm}%
\def\dimHa{{\rm dim_{_{\kern.2mm Ha}}}}
\def\rajou{{}^{}{\Cal B}_{s\,}}
\def\Lis{\Cal L\kern.3mmis\,} 
\def\Linb{\Cal L\lower.7mm\hbox{\kern.1mm\font\=cmmi6\b}}

\def\dualbeta{^{\kern0.4mm\prime}_{\kern-.2mm\raise.95mm\hbox{$_{_\beta}$}}} 
\def\Nbh{\Cal N_{\font\=cmmi6\lower.15mm\hbox{\kern.1mm\bh\kern.15mm}}}
\def\Topma{\roman{{Top_{}}_{\hbox{\font\=cmr6\ma}}\kern.15mm}}
\def\prodc{\prod{_{_{\kern-.3mm\bold c\kern.15mm}}}}
\def\vsprod_#1_#2{\prod\kern-0.3mm{}_{_{\roman{#1}\sp{#2}\,}}} 
\def\vscoprod_#1_#2{\coprod\kern-0.3mm{}_{_{\roman{#1}\sp{#2}\,}}} 
\def\expnota^#1]_#2{\,^{#1\,]{_{}}_{\roman{#2}}}} 
\def\co{\hbox{\font\=cmmi12\c}\kern.15mm\lower.15mm\hbox{$_{\rm o}$}}
\def\LL^#1{L\kern0.15mm\raise.4mm\hbox{$^{#1}$}\kern0.15mm}
\def\lll^#1{\ell\kern0.7mm\raise.2mm\hbox{$^{#1}$}\kern0.15mm}
\def\Cinfty{C\kern.4mm\raise.2mm\hbox{$^\infty$}\kern.15mm}
\def\Cinftyzero{C\raise.2mm\hbox{$^{\,\infty}_{\hbox{\font\=cmr5\0}}$}} 
\def\Cper^#1{C^{\,#1}_{\raise.15mm\hbox{\font\=cmr5\per}}} 
\def\Cperinfty{C\raise.2mm\hbox{$^{\kern.55mm\infty}_{\raise.15mm\hbox{\font\=cmr5\per}}$\kern.15mm}} 
\def\Idc{\roman{{Id_{\kern.3mm}}_{c}\,}} 
\def\Catal_#1{{\raise1.3mm\hbox{\font\=cmr5\c}\alpha_{\kern0.3mm}}_{\boldsymbol{\mathcal#1}}\kern.45mm} 
\def\Catom_#1{{\raise1.3mm\hbox{\font\=cmr5\c}\omega_{\kern0.3mm}}_{\boldsymbol{\mathcal#1}}\kern.45mm} 
\def\baral_#1{\bar\alpha\kern.45mm\lower.2mm\hbox{$_{\boldsymbol{\mathcal#1}}\kern.45mm$}} 
\def\barom_#1{\bar\omega\kern.45mm\lower.2mm\hbox{$_{\boldsymbol{\mathcal#1}}\kern.45mm$}} 
\def\catimes{\kern.95mm\raise.45mm\hbox{\font\=cmbsy6\\char'002}\kern-2.2mm\lower.7mm\hbox{\font\=cmr5\c\kern-.15mm a}\kern1.05mm} 
\def\cstimes{\kern.95mm\raise.45mm\hbox{\font\=cmbsy6\\char'002}\kern-2.1mm\lower.7mm\hbox{\font\=cmr5\c\kern-.15mm s}\kern1.2mm} 
\def\opcat{\,\raise1.7mm\hbox{\font\=cmr5\op}} 

\def\bold#1{{\bf#1}}
\def\roman#1{{\rm#1}}
\def\limu_#1{\lim\kern-5.5mm\lower1.5mm\hbox{$_{#1}\ $}}
\def\oseoy{\raise1.9mm\hbox{\kern.5mm\font\=cmr5\o}\kern-1.7mm y}
\def\O{{}^{}\Cal O}
\def\Univ{\hbox{\font\=cmssbx10\U}{}} 
\def\Pows{\Cal P\kern-.4mm_s\kern.3mm}
\def\lei{      {}_{ {}^{\,\downarrow\text{\hskip-2.1mm}       }  }  \cap       }
\def\lei{\hbox{\kern.45mm$_{^\downarrow}\kern-1.280mm\cap\kern.85mm$}}
\def\Ze{Z\!\!\!Z} 
\def\Zepp{{{{Z\!\!\!Z^{\phantom{l}}}^{{}_{{}^{\!}\!+}}}}}  
\def\Zep{{{{{{{{Z\!\!\!Z_{}}_{}}_{}}_{}}}_{{}^{\!+}}}}}  

\def\inve{\lower.85mm\hbox{$^{^-}$}\kern-.5mm{}^\iota}

\def\fvalue{\hbox{\kern.2mm\font\=cmr10\\char'022\kern-.2mm}} 
\def\ffvalue{\hbox{\kern.2mm\font\=cmr7\\char'022\kern-.2mm}} 
\def\image{\hbox{\font\=cmr10\\char'022\kern-1mm\char'022}} 
\def\iimage{\hbox{\font\=cmr7\\kern.3mm\char'022\kern-.7mm\char'022\kern-.3mm}} 
\def\images{\hbox{\font\=cmr10\\char'022\kern-1mm\char'022\kern-1mm\char'022}} 
\def\weco{\kern.15mm\hbox{\font\=cmtt10\\char'054}\kern.4mm} 
\def\cdotn{\kern-.2mm\cdot\kern-.2mm} 
\def\setminusn{\kern-.2mm\setminus\kern-.2mm} 
\def\timesn{\kern-.2mm\times\kern-.2mm} 
\def\ttimes{\hbox{\kern-.2mm${}\times\kern-2.5mm\lower.8mm\hbox{\font\=cmr5\t}\kern1.8mm$}} 
\def\ttimesn{\hbox{\kern-.2mm${}\times\kern-2.5mm\lower.8mm\hbox{\font\=cmr5\t}\kern1.4mm$}} 
\def\ktimes{\hbox{\kern-.2mm${}\times\kern-2.5mm\lower1mm\hbox{\font\=cmr5\k}\kern1.5mm$}} 
\def\vstimes{\kern.95mm\raise.45mm\hbox{\font\=cmbsy6\\char'002}\kern-2.3mm\lower.9mm\hbox{\font\=cmr5\vs}\kern1.05mm} 

\def\Circ{\kern.9mm\hbox{\font\=cmbsy10\\char'016}\kern.9mm}
\def\cardplus{\hbox{$\kern.77mm+\kern-1.95mm\raise.23mm\hbox{$_{_{\roman c}}$}\kern1.33mm$}}%
\def\ordplus{\hbox{$\kern.78mm+\kern-1.97mm\raise.23mm\hbox{$_{_{\roman o}}$}\kern1.22mm$}}%
\def\svs#1{\sbi{\fiveroman{svs\,}#1}} 

\def\Examplee{{\font\=cmssi10\E\kern.15mmx\kern.15mma\kern.15mmm\kern.14mmp\kern.17mml\kern.15mme}\kern.3mm. }
\def\Examples{{\font\=cmssi10\E\kern.15mmx\kern.15mma\kern.15mmm\kern.14mmp\kern.17mml\kern.15mme\kern.15mms}\kern.3mm. }
\def\Remarkk{{\font\=cmssi10\R\kern.15mme\kern.15mmm\kern.15mma\kern.15mmr\kern.15mmk\kern.15mm. }}
\def\Remarkss{{\font\=cmssi10\R\kern.15mme\kern.15mmm\kern.15mma\kern.15mmr\kern.15mmk\kern.15mms}\kern.3mm. }
\def\N{{I\!\!N}} 
\def\No{{I\!\!N\kern-.54mm\lower.15mm\hbox{$_{\rm o}$}}} 
\def\iNo{I\!\!{N_{}}_{\kern-.22mm{\rm o}}} 
\def\Nopot#1{I\!\!N\kern-.54mm\lower.15mm\hbox{$_{\rm o}$}\kern-.7mm{}^{#1}} 
\def\potNo{^{\kern.37mm I\!\!{N_{}}_{\kern-.22mm{\rm o}}}} 
\def\minus{\kern.2mm\lower1.05mm\hbox{$^-$}}
\def\pplus{\raise.22mm\hbox{\font\=cmr5\\char'053}}
\def\mminus{\raise.18mm\hbox{\font\=cmsy5\\char'000}}
\def\plusinftyy{\raise.18mm\hbox{\font\=cmr5\\char'053}\infty}
\def\minusinftyy{\raise.18mm\hbox{\font\=cmsy5\\char'000}\infty}
\def\plusinfty{\lower1.05mm\hbox{$^+$}\infty}
\def\minusinfty{\lower1.05mm\hbox{$^-$}\infty}
\def\Qe{\hbox{$Q\kern-2.6mm\raise.2mm\hbox{\font\=cmssqi8\I}\kern1.7mm$}}
\def\Re{I\!\!R}
\def\Rep{{{I\!\!R^{\phantom{l}}}^{{}_{{}^{\!}\!+}}}}
\def\Repp{{{{{{{I\!\!R_{}}_{}}_{}}_{}}}_{{}^{\!+}}}}
\def\Ce{{\hbox{$C\kern-2.5mm\raise.2mm\hbox{\font\=cmssqi8\I}\kern1.48mm$}}}
\def\imag{\kern.15mm\lower.6mm\hbox{$^{^*}$}\kern-1.8mm\imath\kern.1mm} 

\def\ebit#1{\kern.1mm\hbox{\font\=cmmib8\#1}\kern.2mm} 
\def\ebiF{\kern.1mm\hbox{\font\=cmmib8\F}\kern.5mm} 
\def\ebiT{\kern.1mm\hbox{\font\=cmmib8\T}\kern.6mm} 
\def\ebiU{\kern.1mm\hbox{\font\=cmmib8\U}\kern.5mm} 

\def\fssi#1{\hbox{\font\=cmssi10\#1}\kern0.15mm} 
\def\smb#1{\hbox{\font\†=cmmi8\†#1\kern.3mm}} 
\def\ssmb#1{\hbox{\font\=cmmi6\#1}} 
\def\eCal#1{\kern.1mm\hbox{\font\†=cmbsy8\†#1\kern.4mm}} 
\def\ecal#1{\kern.1mm\hbox{\font\†=cmsy8\†#1\kern.3mm}} 
\def\ncal#1{\kern.1mm\hbox{\font\†=cmsy9\†#1\kern.3mm}} 
\def\vcal#1{\kern-.1mm\vec{\kern.2mm\hbox{\font\†=cmsy7\†#1}\kern.3mm}} 

\def\id{\kern.3mm\roman{id}\kern.7mm}
\def\idv{\hbox{\font\=cmr10\id}\kern.25mm\lower.8mm\hbox{\font\=cmr7\v}\kern.3mm} 
\def\idm{\hbox{\font\=cmr10\id}\kern.25mm\lower.8mm\hbox{\font\=cmr6\m}\kern.3mm} 
\def\seq#1{\langle#1\rangle}
\def\Seq#1{\big\langle#1\big\rangle}
\def\ymp{{}^{}{\Cal N}_o\,}
\def\SemiNor{\Cal S_{_N}\kern0.15mm}
\def\vecs{\upsilon\kern-0.3mm\lower.15mm\hbox{$_s$}\kern0.2mm} 
\def\vecss{\hbox{\font\=cmitt10\v}\kern-0.1mm\lower.15mm\hbox{$_s$}\kern0.2mm} 
\def\bnull#1{\hbox{\font\=cmssbx10\0}{}_{\font\=cmmi6\lower.15mm\hbox{\kern-.1mm\#1\kern.15mm}}} 
\def\bnulla#1_#2{\hbox{\font\=cmssbx10\0}{}_{\font\=cmmi6\lower.15mm\hbox{\#1\kern-.1mm}}\lower.3mm\hbox{$_{_{#2}}$}} 
\def\bzero#1{\hbox{\font\=cmbx10\0}{}_{\font\=cmmi6\lower.15mm\hbox{\kern-.1mm\#1\kern.15mm}}} 
\def\dom{{{}^{}{\rm dom}\,{}_{{}^{}}}}
\def\domm{\kern0.15mm{\rm dom}{^{\kern.3mm\hbox{\font\=cmr6\2}}}\,}
\def\domr#1{\roman{dom}^{\font\=cmr6\raise.0mm\hbox{\kern.3mm\#1}}}

\def\rng{{}^{}{\rm rng}\,{}_{{}^{}}}
\def\Ccinfty{C\lower.3mm\hbox{$\kern-0.2mm_{\roman c}$}\kern-.85mm\raise.3mm\hbox{$^\infty$}\kern0.15mm} 
\def\CPi#1{C\kern-.2mm\lower.05mm\hbox{$_{_\Pi}$}\kern-1.52mm{}^{#1}}
\def\CinftyPi{C\kern.4mm\raise.3mm\hbox{$^\infty$}\kern-3.35mm_{_\Pi}\kern1.45mm}
\def\CinftyS{\Cinfty\kern-3.9mm_{_{\Cal S}}\kern1.45mm}

\def\RHB#1#2{\raise#1mm\hbox{$#2$}} 
\def\LHB#1#2{\lower#1mm\hbox{$#2$}} 

\def\fiveroman#1{\hbox{\font\=cmr5\#1\kern.1mm}}

\def\sixroman#1{\hbox{\font\=cmr6\#1\kern.1mm}}
\def\eightmath#1{\hbox{\font\=cmmi8\{#1}\kern.1mm}}
\def\eightroman#1{\hbox{\font\=cmr8\{#1}\kern.1mm}}
\def\subtext#1{\raise.2mm\hbox{$_{_{\kern0.15mm\roman{#1}}}$}}
\def\subtexT#1{\raise.2mm\hbox{$_{_{\kern0.15mm\hbox{\font\=cmr5\#1}}}$}}
\def\sNor#1{\kern.25mm\lower.38mm\hbox{$_{#1}$}}
\def\sNorr#1{\kern-.2mm\lower.38mm\hbox{$_{#1}$}}
\def\sNoreset_#1{\kern.13mm\lower.83mm\hbox{\font\=cmmi6\C}\kern.32mm\lower.1mm\hbox{$_{^{\emptyset,#1}}$}}
\def\sbi#1{{_{\kern-0.1mm}}_{#1}} 
\def\ais#1_#2{{}_{\font\=cmmi6\lower.15mm\hbox{\kern-.1mm\#1\kern.15mm}}\lower.3mm\hbox{${_{\kern-0.3mm_{#2}}}$}} %
\def\aais#1_#2{\kern.1mm{}_{\font\=cmmi6\lower.25mm\hbox{\kern-.1mm\#1\kern.15mm}}\lower.4mm\hbox{${_{\kern-0.3mm_{#2}}}$}} %
\def\ai#1{{}_{\font\=cmmi6\lower.15mm\hbox{\kern-.1mm\#1\kern.15mm}}} 
\def\yi#1{^{\font\=cmmi6\raise.0mm\hbox{\kern-.1mm\#1\kern.15mm}}} 
\def\ear#1{{}_{\font\=cmr5\lower.15mm\hbox{\kern.1mm\#1}}} 
\def\ar#1{{}_{\font\=cmr6\lower.15mm\hbox{\kern.1mm\#1}}} 
\def\aar#1{_{\font\=cmr6\lower.15mm\hbox{\kern.1mm\#1}}} 
\def\yr#1{^{\font\=cmr6\raise.0mm\hbox{\kern.3mm\#1}}} 
\def\yrai^#1_#2{^{\kern.4mm\hbox{\font\=cmr6\{#1}}}_{\kern.2mm{#2}}}
\def\upparentes#1{^{\kern.2mm\raise.2mm\hbox{\font\=cmr6\\char'050}\kern.1mm{#1}\kern.1mm\raise.2mm\hbox{\font\=cmr6\\char'051}}} 
\def\lupar{\kern.2mm\lower1mm\hbox{$^{^(}$}} 
\def\rupar{\lower1mm\hbox{$^{^)}$}\kern-.15mm} 
\def\yyi#1{^{\font\=cmmi6\lower.6mm\hbox{\kern-.25mm\#1\kern-.05mm}}} 
\def\yyr#1{^{\font\=cmr6\lower.45mm\hbox{\kern-.25mm\#1\kern-.15mm}}} 
\def\yplus{\lower1mm\hbox{$^{^+}$}} 
\def\yminus{\lower1mm\hbox{$^{^-}$}} 
\def\aminus{{\kern.15mm\raise.3mm\hbox{$_{_-}$}\kern-.1mm}}%
\def\yvee{\LHB{.9}{^{^{\,\vee}}}\kern-.3mm} 
\def\ywed{\LHB{.9}{^{^{\,\wedge}}}\kern-.3mm} 
\def\adot{\kern.2mm\hbox{\font\=cmb10\\char'056}}%
\def\ydot{\kern.2mm\raise1.9mm\hbox{\font\=cmb10\\char'056}}
\def\yydot{\kern.2mm\raise1.35mm\hbox{\font\=cmb7\\char'056}\kern.2mm}
\def\yydott{\kern.2mm\raise1.35mm\hbox{\font\=cmb6\\char'056}\kern.2mm}
\def\ClT{{\rm Cl}\kern.25mm\lower.4mm\hbox{$_{\Cal T}$}\kern0.2mm} 
\def\IntT{\sp{\rm Int}\kern.2mm\lower.4mm\hbox{$_{\Cal T}$}\kern0.2mm} 
\def\Cl_taurd#1{\roman{Cl_{}}_{\kern0.37mm\hbox{\font\=cmmi8\\char'034}\kern-0.15mm{_{}}_{rd}\kern0.2mm#1\,}}
\def\Int_taurd#1{\roman{Int_{}}_{\kern0.37mm\hbox{\font\=cmmi8\\char'034}\kern-0.15mm{_{}}_{\roman{rd}}\kern0.2mm#1\,}}
\def\inc{\subseteq}
\def\iinc{\supseteq}
\def\exi#1{\exists\,#1\kern.2mm\,;}
\def\all#1{\forall\,#1\kern.2mm\,;}
\def\imply{\Rightarrow}

\def\spp{\kern0.07mm} 
\def\sp{\kern0.15mm} 
\def\ssp{\kern0.37mm} 
\def\snn{\kern-0.2mm} 
\def\sn{\kern-0.3mm} 
\def\ssn{\kern-0.63mm} 
\def\biggerlineskip#1 {\linebreak\nopagebreak\vskip-4.2mm\vskip.#1mm\nopagebreak\noindent}%
\def\Biggerlineskip#1 {\linebreak\nopagebreak\vskip-4.2mm\vskip#1pt\nopagebreak\noindent}%
\def\nKP#1{$\null$\kern#1mm}
\def\KP#1{\kern#1mm} 
\def\KN#1{\kern-#1mm} 
\def\nhskip#1mm{$\null$\kern#1mm}

\def\mhyppy#1{\null\kern#1mm}

\def\text#1{\hbox{\rm#1}}

\def\VBOX/#1/#2/HEREend{\vbox{#2\vskip-#1mm}\vfill\null\eject}
\def\œ$#1${\hbox{$#1$}} 
\def\"{\"a} \def\"{\"o}
\def\q#1{``\kern0.37mm#1\kern0.37mm"}
\def\newProCla#1\par#2\par{\vskip1.7mm\noindent\bf#1\it#2\vskip1.7mm}
\def\Prooff{{\font\=cmssi10\P\kern.37mmr\kern.37mmo\kern.37mmo\kern.37mmf\kern.37mm. }\rm}

\def\QED{\hfill\hbox{$\ \sqcap\kern-2.45mm\sqcup$}}

\def\noin{\noindent}
\def\Newline{\kern-10mm\newline}
\font\rp=cmr8
\def\eps{\varepsilon}
\def\leu{\raise1.5mm\hbox{\font\=cmmi5\\char'074}\kern.2mm}%
\def\riu{\kern.2mm\raise1.5mm\hbox{\font\=cmmi5\\char'076}}%
\def\Symbol#1Ï{\kern.35mm\hbox{\font\=cmr10\\char'047}\kern.2mm#1\kern.35mm\hbox{\font\=cmr10\\char'047}}%
\def\Symboo#1Ï{\kern.35mm\text{`}\kern.2mm#1\kern.35mm\hbox{\font\=cmr10\\char'047}}
\def\RunMyHead#1#2#3#4{%
 \headline{\ifnum\pageno=\firstpage\hfil%
           \else{\ifodd\pageno{\rp#3\phantom\folio\hfil#4\hfil\phantom{#3}\folio}%
                 \else{\rp\folio\phantom{#2}\hfil#1\hfil\phantom\folio#2}%
                 \fi}%
           \fi}%
 \footline{\ifnum\pageno=\firstpage\hfil{\rp[\,\folio\,]}\hfil%
           \else\hfil%
           \fi}%
}%
\def\bulgin{\noindent$\bullet$ \ \kern.1mm} 
\def\bulgen{\noindent\kern-1.5mm$\bullet$\kern3.95mm} 
\def\subhead#1\par#2\par{\vskip4mm\smallbreak\null\smallskip\vbox{\noindent\bbf#1\hfill\kern1.5mm#2\hfill\phantom{#1}\vskip2.5mm\nopagebreak}\nopagebreak\noindent}
\def\subheadd#1\par#2\par#3\par{\vskip4mm\smallbreak\null\smallskip\vbox{\noindent\bbf#1\hfill#2\hfill\phantom{#1}\vskip1.5mm\centerline{#3}\vskip2.5mm\nopagebreak}\nopagebreak\noindent}
\def\insubsubhead#1\par{\vskip4mm$\null$\hskip2mm{\font\=cmss10\#1}\vskip2mm\noindent}%
\def\binsubsubhead#1#2\par{\vskip4mm{\bf#1.}\hskip5mm{\font\=cmss10\#2}\nopagebreak\vskip2mm\nopagebreak\noindent}%
\def\wave{\hbox{\font\†=cmsy10\†\hbox{\char'164}\kern-2.35mm\hbox{\char'165}\kern.55mm}}
\def\wavee{\hbox{\font\†=cmsy8\†\hbox{\char'164}\kern-2.0mm\hbox{\char'165}\kern.4mm}} 
\def\barmj{\kern.25mm\bar{\hbox{\font\=cmr10\\char'021}}\kern.4mm}

\def\sigrd{\sigma\kern-.2mm\lower.7mm\hbox{\font\=cmr6\r\font\=cmr5\d}\kern.6mm}
\def\ssigrd{\sigma\kern-.2mm\lower.7mm\hbox{\font\=cmr6\r\font\=cmr5\d}\kern-1.7mm\raise1.25mm\hbox{\font\=cmr6\2}\kern1mm}
\def\sssigrd{\sigma\kern-.2mm\lower.7mm\hbox{\font\=cmr6\r\font\=cmr5\d}\kern-1.7mm\raise1.25mm\hbox{\font\=cmr6\3}\kern1mm}
\def\sigrdu^#1{\sigma\kern-.2mm\lower.7mm\hbox{\font\=cmr6\r\font\=cmr5\d}\kern-1.7mm\raise1.25mm\hbox{\font\=cmr6\#1}\kern1mm}
\def\taurd{\tau\kern-.4mm\lower.7mm\hbox{\font\=cmr6\r\font\=cmr5\d}\kern.6mm}
\def\tsigrd{\tau\sigma\kern-.2mm\lower.7mm\hbox{\font\=cmr6\r\font\=cmr5\d}\kern.6mm}
\def\tauRe{\tau{_{\kern-0.6mm}}_{\hbox{\font\=cmmi5\I\!\!R}}} 
\def\tauR#1{\tau_{_{I\!\!R}}\kern-1.5mm^{#1}}
\def\RN{I\!\!R\kern.3mm^{\hbox{\font\=cmmi6\N}}} 
\def\QTN{Q\kern.1mm_{\lower.2mm\hbox{\font\=cmmi6\T}}^{\kern.2mm\hbox{\font\=cmmi6\N}}} 

\def\leLCS-{{\le}{}_{_{{\rm LCS}}}\text{\sp-\sp}}
\def\Centerline#1\par#2\par#3{\noindent#1\phantom{#3}\hfill#2\hfill\phantom{#1}#3}

\def\rGatDer_#1#2{\raise1.3mm\hbox{\font\=cmr5\r\kern.2mm Gat}\hbox{\font\=cmssi10\D\kern.3mm}\lower.3mm\hbox{$_{#1\kern.2mm #2}$}\kern.5mm} 
\def\nLinb_#1{\mathcal L\kern.3mm\lower.3mm\hbox{$_{\hbox{\font\=cmr6\b}\kern.3mm #1}$}\kern.3mm} 
\def\svs#1{\sbi{\fiveroman{svs\,}#1}}

\begin{document}

\title[$\text{\sc Inverse and implicit function theorems}$]%
            {On Yamamuro's inverse and implicit function \vskip1mm
                 theorems in terms of calibrations}

\author[S. Hiltunen]{Seppo\ I\. Hiltunen}
\address{Helsinki University of Technology                             \vskip0mm$\hspace{2mm}$
           Institute of Mathematics, U311                              \vskip0mm$\hspace{2mm}$
           P.O.\ Box 1100                                              \vskip0mm$\hspace{2mm}$
           FIN-02015 HUT\vskip0mm
         FINLAND}
\email{shiltune\,@\,cc.hut.fi}



\begin{abstract}

For the Fr\'echet space \œ$E=C^{\,\infty}(\sp\roman S\sp^{1.\sn})$ and for a
smooth \œ$\varphi:\Re\to\Re\,$, we prove that the associated map \œ$E\to E$
given by \œ$x\mapsto\varphi\circ x$ satisfies the continuous B$\Gamma\,
$--\,differen- tiability condition in Yamamuro's inverse function theorem only
if $\varphi$ is affine. Via more complicated examples, we also generally
discuss the importance of testing the applicability of proposed inverse and
implicit function theorems by this kind of simple maps.       \end{abstract}

\maketitle


\noin In \cite[p.\ 3]{SeBGN}\ssp, we mentioned that in \cite{Yam79} quite
special differentiabilities are designed hoping to get inverse and implicit
function theorems (see \cite[5.2, 5.3, p.\ 45]{Yam79}\sp) applicable to maps
of Fr\'echet function spaces. Our Theorem \ref{main result} together with the
examples and remarks below indicates this hope to be overoptimistic. Before
getting into the proof of Theorem \ref{main result} in B below, we discuss the
general relevance of this kind of results to refuting applicability of
proposed inverse or implicit function theorems, shortly IFTs\ssp. For the
notations neither immediately guessed by the reader nor explained below, we
refer to \cite[pp.\ 4\,--\,6]{SeBGN} and \cite[pp.\ 4\,--\,9]{Hic}\ssp.


  \binsubsubhead A{Introductory and motivating considerations}\label{subsec A}

We first note that inverse (inFT) and implicit (imFT) function theorems are to
some extent complementary parts of more general IFT type results. Assume that
we are given a class $\ssp\Cal C\aar 1$ of differentiable maps of a certain
order, and loosely say that a function $\sp f\sp$ is {\it regular\ssp} if{}f
we have \œ$\sp(E\ssp,F\sp,f\ssp)\in\Cal C\aar 1$ for some implicitly
understood spaces $\sp E\ssp,F\sp$. Note that $\sp u\fvalue z\sp$ is the
function value of $\sp u\sp$ at $\sp z\ssp$, which conventionally is denoted
by \q{\sp$u\ssp(z)\sp$}\sp. We also have $\ssp f\sp\inve\image\snn B =
\{\,x:\exi{y\in B}\,(\sp x\ssp,y\sp)\in f\sp\,\}\,$.

Now, first suppose that we have an imFT for functions \œ$\sp
f\inc A\timesn B\timesn B\sp$ where $\sp A\sp$ and $\sp B\sp$ are subsets of
structured (e.g.\ topological\sp/\ssp locally convex\sp/\sp normed\sp) vector
spaces $\sp E\sp$ and $\ssp F\sp$, respectively. The imFT says that under
suitable conditions for given \œ$\sp(\sp x\ar 0\ssp,y\ar 0\ssp,b\sp)\in f\sp$,
there is a regular $\sp g\sp$ with \œ$\sp(\sp x\ar 0\ssp,y\ar 0) \in g \inc
f\sp\inve\image\sn\{\sp b\sp\}\,$, and hence we have \œ$\sp
f\sp\fvalue(\sp x\ssp,g\fvalue x\sp)=b\sp$ for all \œ$\sp x\in\dom g\,$.
Suppose further that we have a function\Newline $\sp h\sp$ with \œ$\sp
(\ssp y\ar 0\ssp,x\ar 0)\in h\inc B^{\sp\times 2.}=B\timesn B\sp$ and that we
would be pleased with getting a regular $\sp g\sp$ with \œ$\sp
(\sp x\ar 0\ssp,y\ar 0)\in g\inc h\inve\sp$, a local right inverse to $\sp
h\ssp$. If suitable conditions are satisfied, in the imFT we may take \œ$\sp
E=F\sp$ and \œ$\sp b=\bnull E$ and $\sp f$ given by the pre- scription $\,
(\sp x\ssp,y\sp)\mapsto x-h\fvalue y\ssp$, and get $\sp g\sp$ as required.

Conversely, suppose  we have an inFT and that for a given function $\sp f\sp
$ we want to establish an implicit function $\sp g\sp$ as above. Then (under
suitable conditions) with \œ$\sp B\ar 1=A\times B\sp$ we may take \œ$\sp
h=[\,\sp\roman{pr}\ar 1\ssp,\sp f\sp\,]\subtext f\inc B\ar 1{\sn^{\times 2.}}$
given by \œ$\sp(\sp x\ssp,y\sp)\mapsto(\sp x\ssp,f\sp\fvalue(\sp x\ssp,y\sp))$ \linebreak
in the inFT obtaining \œ$\sp g\ar 1\inc h\inve\sp$, and finally get \œ$\ssp
  g = \roman{pr}\ar 2\snn\circ g\ar 1\snn\circ
      [\,\id,\sp\Univ\times\sn\{\sp b\sp\}\sp\,]\subtext f\sp$ given by the
prescription $\,
   x\mapsto\roman{pr}\ar 2\snn\circ g\ar 1\KN1\fvalue(\sp x\ssp,b\sp)\,$.

\begin{remarks}

As formulated above, from a local imFT one can only get a local right inFT\sp.
However, usually one has such a topological situation that existence of some $\sp
W$ is guaranteed so that in the imFT we may take \œ$\ssp
g = f\sp\inve\image\sn\{\sp b\sp\}\cap\sp W$ with $\sp W$ a neighbourhood of $\sp
(\sp x\ar 0\ssp,y\ar 0)\sp$ in the product topology, cf.\
\cite[Sec.\ 4, Theorems 1\ssp, 5\ssp, pp.\ 19, 20]{Hic}\ssp. From this
stronger formulation of an imFT we get a \q{two\ssp-\ssp sided} local inFT
where existence of $\ssp V$ is guaranteed such that \œ$\sp y\ar 0\in V$ and $\,
h\,|\,V\sp$ injective and regular with also $\sp(\sp h\,|\,V\ssp)\inve$
regular. See, e.g.\ \cite[Corollary 4.6, p.\ 22]{Hic}\ssp.

Further, it should be noted that proving an imFT directly may permit one to
get a result more general than one would get from a previous inFT via our
observations above, and similarly with the roles reversed. For example,
considering the classical Banach or normed space calculus, if one first gets
the inFT\sp, one must consider maps $\sp(E\ssp,E\ssp,h\sp)\sp$ where $\sp E\sp
$ is a {\it complete\ssp} normed space. From this inFT one can only get an
imFT for functions \œ$\sp f\inc A\timesn B\timesn B\sp$ where also the
\q{parameter} set $A\sp$ lies in a Banach space. As for the converse
situation, from our imFT in \cite[p.\ 235]{HiSM} we can directly get only the
classical inFT\sp.                                             \end{remarks}

To make the preceding more concrete, we next consider some examples. For the
definition of the Gateaux derivative function $\sp\rGatDer_EFf\sp$ of a map $\sp
(E\ssp,F\sp,f\ssp)\sp$ see \biggerlineskip3 Definitions \ref{Def Gat} below.
Also note that \œ$\sp 1\adot=\{\emptyset\}\sp$ and \œ$\ssp 2\adot =
\{\ssp\emptyset\,,1\adot\ssp\}\sp$ and \œ$\sp\{\,i\sp\ydot\sn:i\in\No\ssp\} = $
\œ$\Zep\inc\Re\ssp$ and \œ$\sp\sp\{\,n\adot\sn:n\in\Zep\ssp\}=\No$ with \œ$\sp
(\sp n\spp\adot)\ydot=n\sp$ and \œ$\sp(\sp i\sp\ydot\ssp)\adot=i\sp$ and \œ$\ssp
i\ssp\yplus =i+1\adot$ $=i\ssp\cup\sp\{\ssp i\ssp\}\sp$ and $\ssp
i\ssp\yplus\sn\yplus = (\sp i\ssp\yplus)\spp\yplus$ for $\sp i\in\No$ and $\ssp
n\in\Zep\,$.

\begin{example}\label{Exa1}

Existence and \q{regular} dependence on parameters (including initial\sp/\ssp
boundary values and the \q{equations} themselves) of solutions to partial
differential equations can be obtained by using IFTs\ssp. To get a simple
particular case of this general and vague scheme, we consider a partial
problem of the more general one already solved in
\cite[Section 5, Example 7, Theorem 8, pp.\ 30\,--\,31]{Hic}\ssp.

Namely, let \œ$\sp I=[\,0\,,1\,]\sp$ and \œ$\sp Q=I\timesn\Re\sp\,$, and let
the fixed smooth \œ$\sp\varphi:Q\timesn\Re\to\Re$ be such that \œ$\sp
\varphi\fvalue(\sp t\ssp,\eta+1\ssp,\xi\sp)=
\varphi\fvalue(\sp t\ssp,\eta  \ssp,\xi\sp)\sp$ holds for all \œ$\sp
(\sp t\ssp,\eta\ssp,\xi\sp)\in Q\timesn\Re\sp\,$. Letting $\sp E =
\Cperinfty(\Re\ssp)\,$, see the few lines just before Lemma \ref{Le1} below,
and with \vskip.3mm $\null\hfill
      S = \{\,x:\all{t\in I\ssp,\sp\eta\in\Re\sp}\,
              x\fvalue(\sp t\ssp,\eta+1\sp)=
              x\fvalue(\sp t\ssp,\eta  \sp)\ssp\}\null\hfill$ also letting \vskip.3mm

\noin $F\aar 0=C^{\,1.}(Q)_{\sp/\ssp S}\,$ and $\,F=\Cinfty(Q)_{\sp/\ssp S}\,
$, \,assume that $\,x\ar 0\in\vecs E\sp$ and $\,y\ar 0\in\vecs F\sp$ with \vskip.3mm \centerline{$
  \partial\ar 1\ssp y\ar 0+\partial\ar 2\,y\ar 0
= \varphi\circ[\,\id,y\ar 0\,]\subtext f\sp$ and $\,
  y\ar 0\sp(\sp 0\,,\sp\cdot\ssp) = x\ar 0\,$.} \vskip.3mm

\noin In other words, we have a simple nonlinear partial differential equation
on a compact cylinder with boundary values $\sp x\ar 0$ specified on one of
the two components. We are interested to know (?)\ whether there is an open
neighbourhood $\ssp U$ of $\sp x\ar 0$ in the space $\sp E\sp$ such that for
every $\sp x\in U\sp$ there is a unique solution $\sp y\in\vecs F\aar 0\,$,
and that in fact $\sp y\in\vecs F\sp$ and also this correspondence $\sp
x\mapsto y\sp$ defines a smooth map $\ssp g:E\to F\sp$.

From \cite[Theorem 5.8]{Hic} it follows that the answer to (?)\ is
affirmative. However, for the purposes of the note at hand, we sketch another
approach via an inFT as follows. As we saw in
\cite[Example 5.5, pp.\ 26\,--\,27]{Hic}\ssp, we can write the equation in
question for the unknown $\sp y\sp$ with \q{initial} values $\sp x\sp$ as \œ$\,
y=\roman S\,x+\roman I\,(\sp\varphi\circ[\,\id,y\,]\subtext f)\,$, \linebreak
where \œ$\, \roman S\,x = \seq{\,x\fvalue(\sp\eta-t\sp) : \zeta = (\sp t\ssp,
\eta\sp)\in Q\,}\,$ and \œ$\, \roman I\,v =
  \Seq{\sp\int_{\,0}^{\,{\sigma_{\sn}}_{\roman{r\KN{.1}d}}\sp\zeta}
     v\circ{?\sbi{\eightmath\tauu}}\KN{.5}\fvalue\zeta
        \,\ssp\roman d\,\tau:\zeta\in Q\,}\,$ \linebreak with $\,
?\sbi{\eightmath\tauu} =
   \seq{\,(\sp\tau\sp,\eta-t+\tau\sp):\zeta=(\sp t\ssp,\eta\sp)\in Q\,}\,$.

Now, if with \œ$\ssp G=E\sp\sqcap F\sp$ and \œ$\ssp f =
\seq{\,\ssp y-\roman S\,x-\roman I\,(\sp\varphi\circ[\,\id,y\,]\subtext f):z =
(\sp x\ssp,y\sp)\in\vecs G\sp\,}$\Newline and \œ$\,
h=[\,\sp\roman{pr}\ar 1\ssp,\sp f\sp\,]\subtext f\sp$ and \œ$\,\tilde h =
(\sp G\sp,G\sp,h\sp)\,$, we have an inFT applicable to the map $\ssp\tilde h\,
$, we can get a local {\it solution map\ssp} $(E\ssp,F\sp,g\sp)$ as explained
above. However, this does not give uniqueness of the solution. Moreover, we
shall see in Remarks \ref{Rems2}\ssp(b) below that at least the inFT
\cite[5.2, p.\ 45]{Yam79} {\it is not\ssp} applicable in the generic case
where the equation is not linear near $\ssp y\ar 0\,$, i.e.\ when no \œ$\ssp
\eps\in\Rep$ exists with \œ$\,
\partial_{\ssp\sixroman 3}^{\,\ssp 2.}\sp\varphi\image N\inc\{0\}$ \linebreak
for $\,N=Q\timesn\Re\ssp\cap\{\ssp(\sp\zeta\ssp,\xi\sp) : \exi{\xi\ar 1}\,
(\sp\zeta\ssp,\xi\ar 1)\in y\ar 0\text{ and }\ssp|\,\xi-\xi\ar 1\ssp|<\eps\,\}
\,$.                                                           \end{example}

\begin{example}\label{Exa3}

Let \œ$\sp\eps\ar 0\in\Rep$ and \œ$\sp\smb N\in\N\ssp$, and let \œ$\sp Q\ar 0=
{\ssp]}\minus\eps\ar 0\ssp,\eps\ar 0\,{[\ssp}\times\bar\Omega\sp$ where $\sp
\bar\Omega\sp$ is the closure of a bounded open $\sp\Omega\sp$ contained in $
\Re\sp\,\yi N\sp$. Let \œ$A\inc 1\adot\timesn\Nopot{\ssmb N}\sp$ be finite,
and let \œ$\sp\varphi\subtext{old}\sn:Q\ar 0\timesn\Re\sp\,^A\to\Re\,$ be
smooth. For \œ$\sp\bar\alpha=(\sp i\ssp,\alpha\sp)\in\No\timesn\Nopot{\ssmb N}\sp
$ and for a smooth function \œ$\sp y:\bar\Omega\ssp(\smb T\sp) =
[\,0\,,\smb T\sp\,]\times\bar\Omega\to\Re\ssp$ with \œ$\ssp\smb T\in\Rep\sp$,
now the iterated partial derivative \œ$\,\partial\,^{\bar\alpha}\ssp y =
\partial_{\sp\fiveroman T}^{\,i}\sp
\partial_{\sp\fiveroman S}^{\,\sp\alpha}\ssp y\sp$ is defined in an \q{obvious}
manner. Letting the jet \œ$\,\roman J\,y:\bar\Omega\ssp(\smb T\sp) \to
\Re\sp\,^A\sp$ be defined by \œ$\sp \zeta = (\sp t\ssp,\eta\sp)\mapsto
\seq{\,\sp\partial\,^{\bar\alpha}\sp y\fvalue\zeta:\bar\alpha\in A\sp\,}\,$,
we may consi- der the partial differential equation \,oldE\œ$\,
(\ssp y\subtext{ini}\ssp,y\subtext{old}) :
\partial_{\sp\fiveroman T}^{\ssp 1.}\sp y\subtext{old} = \varphi\subtext{old}\sn
\circ[\,\id,\roman J\,y\subtext{old}\,]\subtext f\sp$ \linebreak with initial
condition $\,y\subtext{old}\sp(\sp 0\,,\sp\cdot\,) = y\subtext{ini}\,$.

Assume further that we are given suitable boundary conditions in the form of a
linear subspace $\sp S\ar 0$ in the Fr\'echet space $\sp\Cinfty(\bar\Omega)\,
$. Putting \œ$\sp E\subtext{ini}=\Cinfty(\bar\Omega)_{\sp/\ssp S_0}$ and \œ$\,
\roman F\ssp(\smb T\sp)=\Cinfty(\sp\bar\Omega\ssp(\smb T\sp))_{\sp/\ssp S}\,$,
where \œ$\ssp S=\Re\ ^{\bar\Omega\ssp(\ssmb T\ssp)}\cap\{\,y:\all{t\in[\,0\,,
\sp\smb T\sp\,]\sp}\,y\ssp(\sp t\ssp,\sp\cdot\,)\in S\ar 0\ssp\}\,$, \biggerlineskip2
for a given \œ$\ssp y\subtext{ini\sp 0}\in\vecs E\subtext{ini}$ suppose that
we are interested to know (?)\ whether some $\sp\eps\sp$ with \œ$\sp
0<\eps<\eps\ar 0$ and an open neighbourhood $\sp U$ of $\ssp
y\subtext{ini\sp 0}$ in $\sp E\subtext{ini}$ and a smooth function $\,
g:E\subtext{ini}\iinc U\to\roman F\ssp(\eps)\sp$ exist with $\,
g\inc\{\ssp(\ssp y\subtext{ini}\ssp,y\subtext{old}) :
             \roman{oldE}\,(\ssp y\subtext{ini}\ssp,y\subtext{old})\ssp\}\,$.

In other words, we want to know whether our initial\ssp-\sp boundary value
problem is in a certain sense \q{locally well\ssp-\sp posed}. We approach the
problem via an inFT as follows. Let \œ$\sp E=\bosy R\sp\sqcap E\subtext{ini}$
and \œ$\sp F=\roman F\ssp(1)\sp$ and \œ$\ssp G=E\sp\sqcap F\sp$. Also let \œ$\,
\roman S\,y\subtext{ini} = $ \œ$
    y\subtext{ini}\sn\circ\roman{pr}\ar 2\,|\,\bar\Omega\ssp(1)\sp$ and \œ$\,
\roman J\ar 0\,y =
\Seq{\sp\int_{\,0}^{\,t}y\fvalue(\sp\tau\sp,\eta\sp)\,\ssp\roman d\,\tau :
\zeta=(\sp t\ssp,\eta\sp)\in\bar\Omega\ssp(1)\,}\,$. Considering the \Biggerlineskip1
family \œ$\,
\varphi=\seq{\,\seq{\,\ssp\smb T\sp\,\varphi\subtext{old}\KN{.6}\fvalue
(\ssp\smb T\sp\,t\ssp,\eta\ssp,\xi\sp) : \smb P=(\sp t\ssp,\eta\ssp,\xi\sp)\in
  \bar\Omega\ssp(1)\timesn\Re\sp\,^A\sp\,} :
    \smb T\in{\ssp]}\minus\eps\ar 0\ssp,\eps\ar 0\,{[}\ }\,$, if \linebreak
with \œ$\sp 0<\eps<\eps\ar 0$ we let \œ$\,
\Iota=\seq{\,\seq{\,\sp y\fvalue(\sp\eps^{\sp\mminus 1}\sp t\ssp,\eta\sp) :
\eps : \zeta=(\sp t\ssp,\eta\sp)\in\bar\Omega\ssp(\eps)\,}:y\in\vecs F\sp\,}\,
$, then $\Iota\sp$ is a linear homeomorphism \œ$\sp F\to\roman F\ssp(\eps)\sp$
such that for \œ$\ssp y\in\vecs F\sp$ satisfying the equa- tion $\,
 y = \roman S\,y\subtext{ini}\sn +
  \roman J\ar 0\ssp(\sp\varphi\fvalue\eps\circ[\,\id,\roman J\,y\,]\subtext f)$
also $\,\roman{oldE}\,(\ssp y\subtext{ini}\ssp,\sp\Iota\fvalue y\sp)\,$ holds. \vskip.3mm

Consequently, we get (?)\ affirmatively answered if with \vskip.3mm $\nKP{4.6}
\smb W\sn\ar 0 = (\sp 0\,,y\subtext{ini\sp 0}\ssp,\bnull F\ssp,z\ar 0)\,$, \,
  where $\,z\ar 0 = (\sp x\ar 0\ssp,y\ar 0) = (\sp 0\,,
    y\subtext{ini\sp 0}\ssp,\roman S\,y\subtext{ini\sp 0})\,$, \hfill and \KP{11} \vskip.5mm

$\nKP7 h = \seq{\,(\sp x\ssp,y - \roman S\,y\subtext{ini}\sn - \roman J\ar 0\ssp
      (\sp\varphi\fvalue\ssp\smb T\sp\circ[\,\id,\roman J\,y\,]\subtext f)) : $ \newline $\null\hfill
  z = (\sp x\ssp,y\sp) = (\ssp\smb T\sp,y\subtext{ini}\ssp,y\sp)\in\vecs G\sp$
      and $\,|\,\smb T\,| < \eps\ar 0\,}$ \KP{10} \vskip.5mm

\noin we show existence of $\ssp h\ar 1$ with \œ$\ssp\smb W\sn\ar 0\in h\ar 1
\inc h\inve$ and $\sp(\sp G\sp,G\sp,h\ar 1)\sp$ smooth. Namely, then there are
\œ$\sp\eps\in\Rep$ and an open neighbourhood $\sp U$ of $\ssp
y\subtext{ini\sp 0}$ in $\sp E\subtext{ini}$ such that we have \œ$
[\minus\eps\ssp,\eps\,]\times U\times\snn\{\ssp\bnull F\}\inc\dom h\ar 1\ssp$,
and we may take \œ$\,g=\Iota\circ g\ar 1$ where $\ssp g\ar 1$ is defined by
the prescription $\,U \owns  y\subtext{ini} \mapsto \roman{pr}\ar 2\snn\circ
    h\ar 1\KN1\fvalue(\sp\eps\ssp,y\subtext{ini}\ssp,\bnull F\spp)\,$.

We shall see in Remarks \ref{Rems2}\ssp(c) below that \cite[5.2, p.\ 45]{Yam79}
is of no use here pro- \Biggerlineskip1 vided the boundary conditions satisfy \œ$\ssp
\vecs C\,(\bar\Omega)\cap\{\,v:v\,|\,\Omega\in\vecs\Cal D\ssp(\Omega)\,\} \inc
S\ar 0\,$, and the \q{equation} or the pair \œ$\ssp \eCal E =
(\sp\varphi\subtext{old}\ssp,\sp y\subtext{ini\sp 0})\sp$ is {\it initially
strictly nonlinear\ssp}, meaning that some \œ$\sp\eta\ar 0\in\Omega\ssp$
exists such that for every $\sp\eps\sp$ with \œ$\ssp 0<\eps\le\eps\ar 0$ there
are \œ$\,\bar\alpha\,,\sp\bar\alpha\ar 1\in A\,$ and \œ$\,t > 0\sp$ and \œ$\,
\eta\in\Omega\ssp$ and \œ$\,\xi\in\Re\sp\,^A$ with \œ$\,
t + |\,\eta - \eta\ar 0\ssp|\subtexT{\char'006}\snn + |\,\xi - \roman{J\,S\,}
y\subtext{ini\sp 0}\KN{.6}\fvalue(\sp t\ssp,\eta\sp)\,|\subtexT{\char'006}
< \eps$ \Biggerlineskip1 and $\,
   \partial_{\,\bar\alpha}\,\partial_{\,\bar\alpha_1}\sp
    \varphi\subtext{old}\KN{.5}\fvalue(\sp t\ssp,\eta\ssp,\xi\sp) \not= 0\,$,
\,generally letting $\,|\ssp\zeta\ssp|\subtexT{\char'006} =
      \sum\,\seq{\,\sp|\,\zeta\sp\fvalue i\,|:i\in\dom\zeta\sp\,}\,$.

Note that according to the preceding definition from $\sp\eCal E\sp$ not being
initially strict- ly nonlinear it follows existence of $\ssp N\aar 2$ with the
property \œ$\,
\roman{J\,S\,}y\subtext{ini\sp 0}\ssp|\,(\spp\{0\}\sn\times\Omega\ssp) \inc
N\aar 2 \in$ \œ$\taurd(\bosy R\sp\sqcap\bosy R\sp\expnota^\ssmb N]_{tvs}
                              \sn\sqcap\bosy R\sp\expnota^      A]_{tvs})
      \lei (\Repp\!\times\Omega\times\sn\Re\sp\,^A\sp)\,$, and also such that
for every fixed \œ$\zeta\in\dom N\aar 2$ there is an affine \œ$\sp
\Alf:\bold R\expnota^A]_{vs}\to\bold R\ssp$ with \œ$\,
\varphi\subtext{old}\ssp|\,N\aar 2\ssp(\sp\zeta\,,\sp\cdot\,)\inc\Alf\,$. That
is, then the equation is in a sense locally linear near the initial values.
 \end{example}

\begin{example}

Letting \nopagebreak \vskip.3mm $\,
[\ssp x\ssp]\subtext{S1} = \{\,y : x\ssp,y\in D\aar 0\text{ and }\ssp\exi{n\in
  \Ze\sp}\,\all{s\in\Re\sp}\,x\fvalue s=y\fvalue s+n\,\}\,$, \hfill where \vskip.3mm $\nKP{3.9}
D\aar 0 = \vecs\Cinfty(\Re\ssp)\cap\{\,x:0\not\in\rng\roman D\,x\ssp$ and $\,
          \all{s\in\Re\sp}\,x\fvalue(\sp s+1\sp)=x\fvalue s+1\ssp\}\,$, \vskip.3mm

\noin with $\nKP{1.6}
D = \{\,\sp[\ssp x\ssp]\subtext{S1}\ssn : x\in D\aar 0\ssp\}\,$ also putting \vskip.5mm $\nKP{6.3}
\gamma = D^{\sp\times 3.}\sn\cap\{\ssp(\sp\hat x\ssp,\hat y\ssp,\hat z\sp) :
\exi{x\ssp,y\ssp,z}\,(\sp x\ssp,y\ssp,z\sp)\in\hat x\timesn\hat y\timesn\hat z
 \sp\text{ and }\,z=y\circ x\,\}\,$, \vskip.3mm

\noin then $\sp\gamma\sp$ is a group\ssp(\,\ssp operation\sp) on $\sp D\ssp$.
Also let \œ$\Iota=\id\Re\sp$ and \œ$\ssp E = \Cperinfty(\Re\ssp)\sp$ and $
U\aar 1 = \vecs\Cper^{1.}(\Re\ssp)\spp \cap \{\,x : |\,x\fvalue 0\,|<\frac 12\sp$
and $\minus 1\not\in\rng\roman D\,x\,\}\,$ and $\,U=\vecs E\ssp\cap\sp U\aar 1\ssp
$ and \par $\nKP{54} \phi \sp = \sp
\seq{\,\ssp[\,\sp u+\Iota\sp\,]\subtext{S1}\ssn : u\in U\sp\,}\inve\sp$. \vskip.3mm

Considering the smoothness \œ$\sp\Cal S=\CPi{\sp;\,\infty}\sp(\bosy R\ssp)$
corresponding to the differentiability class \œ$\sp\CinftyPi(\bosy R\ssp) =
\Ccinfty(\bosy R\ssp)\,$, by a result similar to
\cite[Ch.\ 3, Sec.\ 1.9, Proposition 18, pp.\ 226\,--\,227]{Bou} and
\cite[Proposition 1.13, pp.\ 354\,--\,355]{Gl-JFA}\ssp, if certain conditions
($*$) hold, there is a unique $\sp M\sp$ with $\sp(M\sp,\gamma\sp)\sp$ a Lie $\sp
\Cal S\,$--\,group such that for some $\ssp V$ with $
[\,\Iota\ssp]\subtext{S1}\sn\in V\sp$ also $\,
M\sp\cup\{\ssp(\sp\phi\,|\,V\spp,E\ssp)\ssp\}\sp$ is an $\,\Cal S\,$--\,atlas.

One of the conditions ($*$) holds if for \œ$\,
       g = \sp U\sp^{\times 3.}\snn\cap\{\ssp(\sp x\ssp,y\ssp,z\sp) :
                      y=z+x\circ(\sp\Iota+z\sp)\ssp\}\sp$ and $\,\tilde g =
(\sp E\sp\sqcap E\ssp,E\ssp,\sp g\ssp)\sp$ we have $\,\tilde g \in
\CinftyPi(\bosy R\ssp)\,$. Letting \vskip.3mm\centerline{$
f = \{\ssp(\sp x\ssp,y\ssp,z\ssp,y-z-x\circ(\sp\Iota+z\sp)) :
         (\sp x\ssp,y\ssp,z\sp)\in U\sp^{\times 3.}\ssp\}\,$,} \vskip.3mm

\noin we see that if we have a suitable imFT applicable to the map \œ$\sp
(\spp E\sp\sqcap E\sp\sqcap E\ssp,E\ssp,f\ssp)\sp$ at the point $\sp
(\sp\bnull E\ssp,\bnull E\ssp,\bnull E\ssp,\bnull E)\,$, we immediately get \œ$\sp
\tilde g\in\CinftyPi(\bosy R\ssp)\,$. Our result
\cite[Theorem 4.3, pp.\ 19\,--\,20]{Hic} is such an imFT\sp, but
\cite[5.3, p.\ 45]{Yam79} is not. Omitting further details, we here only
mention that as the Bp$\ar 2\ssp$--\,extension required in \cite{Hic} one
takes the family $\,\seq{\,(\ssp\roman F\,i\ssp,\roman f\,i\ssp,\idv E\ssp) :
i\in\N\sp\,}\sp$ where $\,\roman F\,i=\Cper^i(\Re\ssp)\sp$ and \vskip.3mm\centerline{$
\roman f\,i = \{\ssp(\sp x\ssp,y\ssp,z\ssp,y-z-x\circ(\sp\Iota+z\sp)) :
                    (\sp x\ssp,y\ssp,z\sp)\in
     U\sn\times U\sn\times(\sp\vecs\roman F\,i\sp\cap U\aar 1)\ssp\}\,$,} \vskip.3mm

\noin and that one gets \œ$\sp
\{\ssp(\spp E\sp\sqcap E\sp\sqcap\spp\roman F\,i\ssp,\roman F\,i\ssp,
\roman f\,i\sp):i\in\N\sp\,\}\inc\CinftyPi(\bosy R\ssp)\sp$ from
\cite[Propositions 0.10, 0.11 p.\ 240]{HiSM} and
\cite[Proposition 3.1, Theorem 3.6, pp.\ 15, 17]{Hic} and the fact \œ$
\{\ssp(\spp G\sp,H\sp,\sp\ell\ssp):G\sp,H\in\LCS(\bosy R\ssp)\text{ and }\ssp
\ell\in\Cal L\,(\spp G\sp,H\ssp)\ssp\}\inc\CinftyPi(\bosy R\ssp)\sp$ by noting
that for $\ssp\roman f\,i\sp$ we have the decomposition \vskip.3mm $\nKP{16}
       (\sp x\ssp,y\ssp,z\sp) \mapsto (\sp\bar x\ssp,y\ssp,z\sp) \mapsto
       (\sp\bar x\circ[\,\id,z\,]\subtext f\ssp,y\ssp,z\sp)$ \par $\nKP{30}
= \ssp (\sp u\ssp,y\ssp,z\sp) \mapsto y-z-u=y-z-x\circ(\sp\Iota+z\sp)\,$, \vskip.3mm

\noin where the prescription \œ$\,x\mapsto\bar x = \{\ssp(\sp s\ssp,t\ssp,
x\fvalue(\sp s+t\sp)):s\ssp,t\in\Re\sp\,\}\,$ defines a continuous linear map
$\,\Cinfty(\Re\ssp)\to\Cinfty(\Re\sp\timesn\Re\ssp)\,$.        \end{example}

The above $\sp(M\sp,\gamma\sp)\sp$ of course is nothing but an interpretation
of the Lie group \ssp Diff$\sp\RHB{.45}{_{_+}}\,\mathbb S^{\,1.}$ of
orientation preserving smooth diffeomorphism of $\ssp\mathbb S^{\,1.}\sp$,
constructed so that one may avoid considering the (quite simple) manifold
structure of $\ssp\mathbb S^{\,1.}\sp$. Considering an arbitrary smooth
finite\ssp-\sp dimensional paracompact (but not necessarily second countable)
manifold $\ssp M\sn\subtext{bas}\ssp$, in a manner similar to the above, we
can use \cite[Theorem 4.3]{Hic} as a tool when constructing the Lie group
\ssp Diff$\sp\,M\sn\subtext{bas}\ssp$. Only the formal details become much
more complicated than in the above simple case. It is our intention to give
them in `\ssp Mapping families, differentiation, and an application to Lie
groups of diffeomorphisms\ssp' although some time will be required for the
completion of this manuscript. In this connection, one should also note
\cite{Mi-Shi} where the same construction problem is treated in a different
manner, however, assuming that the topology of $\ssp M\sn\subtext{bas}$ is
second countable, and still omitting many technical details, although the
presentation there generally is unusually detailed.

\begin{remarks}\label{Rems2}

(a) \  In the proof of Theorem \ref{main result} below, we need a function $u$
   whose $i\ar 0{^{\fiveroman{th}}}$ canonical (semi)\sp norm is small, and
having the absolute value of the $(\sp i\ar 0\yplus){^{\,\fiveroman{th}}}$
derivative large at a given point $\sp s\ar 0\,$, with also $\sp
u\fvalue s\ar 0$ equal to zero. There we can take as $u$ a simple
trigonometric function. For function spaces over more general domains, e.g.\
finite\ssp-\sp dimensional smooth manifolds, we can achieve the same goal by
taking instead as $u$ a suitable scalar multiple of a monomial $m$ multiplied
by a smooth \q{bump} function $b\,$, pulled back by a chart, and extended by
  zero.

More precisely, with \œ$\smb N\in\N$ and \œ$\sp 0<\delta\le 1\ssp$, for \œ$\sp
\alpha\in\Nopot{\ssmb N}$ and \œ$\sp\eta\in\Re\sp\,\yi N$ letting \œ$\ssp
\prod\,\eta=\prod\sp\sbi{i\ssp\in\ssp\ssmb N\,}(\sp\eta\fvalue i\sp)$ and \œ$\ssp
\eta\,^\alpha = \prod\,\seq{\,(\sp\eta\fvalue i\sp)\,^{\alpha\ffvalue i}\sn :
i\in\smb N\,}$ and \œ$\ssp\alpha\ssp! =
\prod\sp\sbi{i\ssp\in\ssp\ssmb N\,}((\sp\alpha\fvalue i\sp)\ssp!\sp)\,$, one \Biggerlineskip1
takes $\, m = \roman m\,^\alpha =
 \seq{\,\eta\,^\alpha\sn:\eta\in\Re\sp\,\yi N\sp\,}\,$ and $\,
 b=\seq{\,\ssp\prod\,(\sp b\ar 0\snn\circ(\sp\delta^{\sp\mminus 1}\sp\eta\sp))
 : \eta\in\Re\sp\,\yi N\sp\,}\,$, where \vskip.5mm \centerline{$
 b\ar 0 = \seq{\,(\sp 1+\exp\,((
            (\sp 2-|\ssp s\ssp|\ssp)^{\,2}\snn-1\sp)^{\sp\mminus 1}\sp
            (\sp 2-|\ssp s\ssp|\ssp)))^{\sp\mminus 1}\sn :
              s\in\Re\sp$ and $\ssp 1<|\ssp s\ssp|<3\sp\,}$} \vskip.2mm $
\mhyppy{51}
  \cup(\ssp[\minus 1\ssp,1\,]\times\sn\{\sp 1\sp\})
  \cup((\Re\setminus{\ssp]}\minus 3\ssp,3\,{\ssp[\sp\,})\timesn\{0\})\,$. \vskip.5mm

Using the Leibniz formula \cite[(2)\ssp, p.\ 101]{Horv}\ssp, with \œ$\sp
\eta\ar 0=\smb N\timesn\{0\}\,$, one sees that then \œ$\ssp
b\cdotn m\fvalue\eta\ar 0=0\,$, unless \œ$\,\alpha=\smb N\timesn 1\adot\ssp$,
and \œ$\,\partial\,^\alpha\sp(\sp b\cdotn m\sp)\fvalue\eta\ar 0 =
                              (\sp\alpha\ssp!\sp)\ydot\sp$, and also derives
existence of \œ$\sp\smb M\in\Rep$ independent of $\sp\delta$ such that for all
\œ$\eta\in\Re\sp\,\yi N$ and for \œ$\kappa\le\alpha$ as functions \œ$
\smb N\to\No$ we have the inequality $\,
    |\,\sp\partial\,^\kappa\sp(\sp b\cdotn m\sp)\fvalue\eta\,|
\le  \delta\sp\,^{|\sp\alpha\sp|\ssp-\ssp|\sp\kappa\sp|}\,\smb M\,$. \par\nKP{-2.5}%
  (b) \ Another aspect of the proof of Theorem \ref{main result} below is that
a contradiction follows from the assumption that a certain {\it nonaffine\ssp}
map $\sp\tilde f\sp$ is \q{\ssp$\roman{C\sbi{B\Gamma}}$} in the sense
\cite[p.\ \nolinebreak 23]{Yam79}\ssp, i.e.\ continuously cb $\,\Gamma\,
$--\,differentiable \œ$\ssp\mu\to\nu$ within $\ssp\bosy E\sp$ in the sense of
our Definition \ref{defi C_{BGamma}} below. More specifically, this
contradiction consists of the formulas \,not \œ$\smb R<\smb R$ and \œ$\,
(*) \ \ \smb R<\smb{A\KP1 D} - \smb M\snn\ar 0\,\smb N\sn\ar 0\le\smb R$ where
\œ$\sp\smb M\snn\ar 0\,,\sp\smb R\in\Repp$ with $\sp\smb M\snn\ar 0$
independent of the varied function $u$ within certain bounds. In ($*$) we have
\q{$\le$} for all $u$ while \q{$<$} only for suitably chosen ones. We have \œ$
\smb A=|\,\varphi\ssp''\fvalue s\ar 0\,|\,$, and hence \œ$\smb A>0$ for
suitably chosen $s\ar 0$ if $\varphi$ is not affine. Further \œ$
\smb D=|\,\sp\roman D\,^{i_0\ssp+\ssp 1.}\ssp u\fvalue s\ar 0\,|\,$ and \œ$
\smb N\sn\ar 0 = \sup\,\{\,|\,\roman D^{\,l}\sp u\fvalue s\,| :
      l\in i\ar 0\yplus\text{ and }\sp s\in\Re\sp\,\}\,$.
The contradiction is obtained by choosing $u$ so that $\smb D$ becomes large
while $\smb N\sn\ar 0$ remains small.

Suppose that instead of $\ssp(E\ssp,E\ssp,f\ssp)$ of Theorem \ref{main result}
as $\tilde f$ we have the $\sp\tilde h\sp$ of Example \ref{Exa1} above. In the
\q{nonlinear} case where \œ$\,
\partial_{\ssp\sixroman 3}^{\,\ssp 2.}\sp\varphi\fvalue\smb P\ar 0\not=0\sp$
for some \œ$\sp\smb P\ar 0=(\sp t\ar 0\ssp,\eta\ar 0\ssp,\xi\ar 0)$ \œ$\in y
\inc Q\timesn\Re\sp\,$, we can prove that $\sp\tilde h\sp$ is not \q{\ssp$
\roman{C\sbi{B\Gamma}}$} by establishing a corresponding inequality \œ$\smb R<
\smb{A\KP1 D} - \smb M\snn\ar 0\,\smb N\sn\ar 0\le\smb R\,$, noting the
following complications. We have $\,\rGatDer_GGh\fvalue(\sp z+w\sp)\fvalue
w\ar 1 = (\sp\bnull E\ssp,v\ar 1\snn - \roman I\,(\sp
\partial\ar 3\ssp\varphi\circ[\,\id,y+v\,]\subtext f\sn\cdot v\ar 1)\,$ when $\,
z\ssp,w\ssp,w\ar 1\in\vecs G$\Newline with \œ$\ssp z=(\sp x\ssp,y\sp)\sp$ and
\œ$\ssp w=(\sp\bnull E\ssp,v\sp)\sp$ and \œ$\ssp
w\ar 1=(\sp\bnull E\ssp,v\ar 1)\,$. We take \œ$\ssp
v\ar 1=Q\timesn\{\sp 1\sp\}\,$, and for the construction of $\ssp v\sp$ we
proceed as follows.

For the space $\sp G\sp$ we consider the \q{canonical} (semi)\ssp norms
\œ$\,\seq{\,\sp\|\ssp w\ssp\|\sbi i\sn:w\in\vecs G\sp\,}\,$ where for $\ssp
w = (\sp u\ssp,v\sp)\in\vecs G\sp$ with $\ssp e = (\sp 1\ssp,1\sp)\sp$
we have \vskip.5mm $\nKP{4.1}
\|\ssp w\ssp\|\sbi i = \sup\,\{\,\sp|\,\sp\roman D\,^{l_1}u\fvalue\eta\sp\,| +
 |\,\sp\partial_{\sp\sixroman 1}^{\,k}\sp\partial_{\sp\sixroman 2}^{\,\sp l}\,
  d\,(\ssp l\ar 0\timesn\{\sp e\})\,v\fvalue\zeta\sp\,| : $ \vskip.3mm $\nKP{28}
    l\ar 1\ssp,\sp k+l+l\ar 0\in i\ssp\yplus$ and $\ssp l\ar 0\in 2\adot$
       and $\ssp\eta\in\Re\sp$ and $\ssp\zeta\in Q\,\}\,$, \vskip1mm

\noin noting that $ \hfill
d\,(\ssp l\ar 0\timesn\{\sp e\})\,v \, = \, \begin{cases} \
d\,(\emptyset)\,v = v & \text{ when }\,l\ar 0=\emptyset\,, \\ \
d\,(\spp\seq{\ssp e\ssp}\spp)\,v = \partial\ar 1\sp v+\partial\ar 2\ssp v & \text{
                               when }\,l\ar 0=1\adot\,,    \end{cases} \hfill\null $ \vskip1mm

\noin and that $\,d\,(\spp\seq{\ssp e\ssp}\spp)\,\roman I\,v=v\ssp$, when $\,
v\in\vecs F\sp$. We take $\sp i=i\ar 0\sn\yplus\sn\yplus$ and \vskip.3mm \centerline{$
v\ssp=\sp\big\{\ssp\big(\ssp t\ssp,\eta\ssp,\sp\delta^{\sp\mminus\frac 12}
  \,b\ar 0\KN1\fvalue(\sp\delta^{\sp\mminus 1}\sp(\sp t-t\ar 0))\cdotn
    (\sp t-t\ar 0)^{\,i_0\ssp+\ssp 1.}\big):t\in I\sp$ and $\ssp
      \eta\in\Re\sp\,\big\}\,$,} \vskip.3mm

\noin where \œ$\sp\delta\in\Rep$ is chosen so that $\ssp
\|\ssp w\ssp\|\sbi{i_0}$ will be small while $\,
   \partial_{\sp\sixroman 1}^{\,\ssp i_0\ssp+\ssp 1.}\,
                   v\fvalue(\sp t\ar 0\ssp,\eta\ar 0)\sp$ becomes large. For
further details, in particular as for the proper order of the various choices,
we refer the reader to the proof of Theorem \ref{main result} below. \par\nKP{-2.5}%
  (c) \ In the situation of Example \ref{Exa3} above, with the provision made
  there at the end, we obtain ($*$) as follows. First, by the nonlinearity
assumption, whatever \œ$ W\sn\in\Nbh(\sp z\ar 0\ssp,\taurd G\ssp)\sp$ we
choose, there always are some \œ$\sp\bar\alpha\,,\sp\bar\alpha\ar 1\in A\sp$
and $\ssp s\ar 0\ssp,\eta\ar 0\ssp,\xi\ar 0\ssp,\smb T\ssp$ and $\ssp y\sp$
with \œ$\ssp 0\le s\ar 0<\smb T\ssp$ and \œ$\,\eta\ar 0\in\Omega\,$, and also
such that for \œ$\ssp\smb P\ar 0 = (\sp s\ar 0\ssp,\eta\ar 0\ssp,\xi\ar 0)\sp$
and \œ$\sp\zeta\ar 0 = (\ssp\smb T^{\,\mminus 1}\sp s\ar 0\ssp,\eta\ar 0)\sp$
and \œ$\ssp z = (\ssp\smb T\sp,\sp y\subtext{ini\sp 0}\ssp,\sp y\sp)\,$, we
have \œ$\ssp z \in W\sp$ and \œ$\ssp
(\ssp\smb T^{\,\mminus 1}\sp s\ar 0\ssp,\eta\ar 0\ssp,\xi\ar 0) \in
\roman J\,y\sp$ and \œ$\,
\partial_{\,\bar\alpha}\,\partial_{\,\bar\alpha_1}\sp\varphi\subtext{old}\KN{.5}
\fvalue\smb P\ar 0\not=0\,$. The contradiction is obtained by considering the
required continuity of $\,\rGatDer_GGh\ssp$ at this point $\ssp z\,$.

We may assume \œ$\ssp|\,\bar\alpha+\bar\alpha\ar 1\ssp|\ssp$ to be the largest
possible. For \œ$\sp \zeta = (\sp t\ssp,\eta\sp) \in
\Re\sp\times\snn\Re\sp\,\yi N\sp$ and \œ$\ssp \bar\kappa=(\sp i\ssp,\kappa\sp)
\in\No\timesn\Nopot{\ssmb N}\sp$ letting \œ$\ssp \zeta\sp\,^{\bar\kappa} =
t\,^i\ssp\eta\,^{\kappa}$ and \œ$\,\bar\kappa\ssp!=i\ssp!\cdotn\kappa\ssp!\,$,
note that if also $\sp \bar\nu \in \No\timesn\Nopot{\ssmb N}\sp$ and $\,\sigma
= \partial\,^{\sp\bar\nu}\ssp\seq{\,(\sp\zeta-\zeta\ar 0)\,^{\bar\kappa}\sn :
   \zeta\in\bar\Omega\ssp(1)\,}\fvalue\zeta\ar 0\,$, \,then \par $\nKP{26}
\sigma = \bar\kappa\ssp !\sp\ydot$ if $\,\bar\kappa=\bar\nu\ssp$, and $\,
\sigma = 0\sp$ otherwise.

By \œ$\sp\eta\ar 0\in\Omega\,$, and by the assumption on the boundary
conditions, we may choose \œ$\sp\delta\ar 1\in\Rep$ so that for \œ$\, v\ar 1 =
\seq{\,\ssp\prod\,(\ssp b\ar 0\snn\circ(\sp\delta\ar 1\KN1^{\mminus 1\,}\taurd
 \snn(\sp\zeta-\zeta\ar 0)))\cdotn(\sp\zeta-\zeta\ar 0)^{\,\bar\alpha_1}\ssn :
    \zeta\in\bar\Omega\ssp(1)\,}$ we have \œ$\sp v\ar 1\in\vecs F\sp$. We take
\œ$\sp i=i\ar 0\sn\yplus\sn\yplus$ and \œ$\ssp
\bar\kappa = (\ssp i\ar 0\sn\yplus\spp,\sp\smb N\timesn 1\adot\spp)\sp$ and \œ$\,
\bar\kappa\ar 1 = (\ssp i\,,\sp\smb N\timesn 1\adot\spp)\,$. We let $\sp
\smb R = j\,\|\,v\ar 1\ssp\|\sbi j\,$, and we choose $\sp\delta\sp$ with $\sp
 0 < \delta \le \delta\ar 1$ so that we get \vskip-.1mm \centerline{$
\smb R + \smb M\snn\ar 0\,\smb N\sn\ar 0 < \smb T\ssp\,|\,\sp
    \partial_{\,\bar\alpha}\,\partial_{\,\bar\alpha_1}\sp\varphi\subtext{old}\KN{.5}
     \fvalue\smb P\ar 0\,|\cdot\delta^{\sp\mminus\frac 12}\ssp
          (\sp\bar\alpha+\bar\kappa\sp)\ssp !\sp\ydot\snn\cdot
            \bar\alpha\ar 1\sp !\sp\ydot\,$,} \vskip.3mm

\noin for certain \œ$\,\smb M\snn\ar 0\,,\sp\smb N\sn\ar 0\in\Rep$ which are
determined once $\ssp\delta\ar 1\ssp,\sp\bar\alpha\,,\sp\bar\alpha\ar 1\ssp,\sp
\bar\kappa\,,\sp\zeta\ar 0\sp$ have been \linebreak fixed. With also \œ$\, v =
  \Seq{\,\delta^{\sp\mminus\frac 12}\ssp\prod\,(\ssp b\ar 0\snn\circ
    (\sp\delta^{\sp\mminus 1\,}\taurd\snn(\sp\zeta-\zeta\ar 0)))\cdotn
                 (\sp\zeta-\zeta\ar 0)^{\,\bar\alpha\ssp+\ssp\bar\kappa}\sn :
\zeta\in\bar\Omega\ssp(1)\,}\,$, we now take $\ssp w=(\sp\bnull E\ssp,v\sp)\sp
$ and $\,w\ar 1=(\sp\bnull E\ssp,v\ar 1)\,$, to obtain \vskip.5mm $\nKP{4.8}
\smb R < \ssp \smb T\ssp\,|\,\sp\partial_{\,\bar\alpha}\,
 \partial_{\,\bar\alpha_1}\sp\varphi\subtext{old}\KN{.5}\fvalue\smb P\ar 0\,|
    \cdot\delta^{\sp\mminus\frac 12}\ssp
       (\sp\bar\alpha+\bar\kappa\sp)\ssp !\sp\ydot\snn\cdot
         \bar\alpha\ar 1\sp !\sp\ydot\snn - \smb M\snn\ar 0\,\smb N\sn\ar 0  $ \vskip.5mm $\nKP{8.35}
 = \spp |\,\sp\partial_{\,\bar\alpha}\,\partial_{\,\bar\alpha_1}(\sp\varphi
    \fvalue\sp\smb T\ssp)\circ[\,\id,\roman J\,(\sp y+v\sp)\,]\subtext f\sn
    \cdot\partial\,\sp^{\bar\alpha\ssp+\ssp\bar\kappa}\ssp v\sn
    \cdot\partial\,\sp^{\bar\alpha_1}\ssp v\ar 1\KN{.6}\fvalue\zeta\ar 0\,|
                                       - \smb M\snn\ar 0\,\smb N\sn\ar 0   $ \vskip.5mm $\nKP{8.35}
 \le \spp |\,\,\partial\,\sp^{\bar\kappa}\ssp(\sp\partial_{\,\bar\alpha_1}(\sp
   \varphi\fvalue\sp\smb T\ssp)\circ[\,\id,\roman J\,(\sp y+v\sp)\,]\subtext f)
    \cdot\partial\,\sp^{\bar\alpha_1}\ssp v\ar 1\KN{.6}\fvalue\zeta\ar 0\,|$ \vskip.8mm $\nKP{8.35}
 = \big|\,\,\partial\,\sp^{\bar\kappa}\ssp\big(\sp
   \sum_{\,\sp\bar\nu\ssp\in\ssp A\,}\partial_{\,\spp\bar\nu\,}(\sp
   \varphi\fvalue\sp\smb T\ssp)\circ[\,\id,\roman J\,(\sp y+v\sp)\,]\subtext f\sn
   \cdot\partial\,\sp^{\bar\nu}\sp v\ar 1\sp\big)\fvalue\zeta\ar 0\,\big|$  \vskip.7mm \centerline{$
 = \big|\,\,\partial\,\sp^{\bar\kappa_1}\sp\big(\ssp
   v\ar 1 + \spp\roman J\ar 0\ssp\big(\sp\sum_{\,\sp\bar\nu\ssp\in\ssp A\,}
   \partial_{\,\spp\bar\nu\,}(\sp\varphi\fvalue\sp\smb T\ssp)\circ
   [\,\id,\roman J\,(\sp y+v\sp)\,]\subtext f\sn\cdot
   \partial\,\sp^{\bar\nu}\sp v\ar 1\sp\big)\big)\fvalue\zeta\ar 0\,\big|$} \vskip.7mm $\nKP{8.35}
 \le \sp |\sn|\sn|\,
       \rGatDer_GGh\fvalue(\sp z+w\sp)\fvalue w\ar 1\,|\sn|\sn|\sp\sbi i
 \le \smb R\,$. \vskip.7mm

Again, to get a proper proof from the preceding pieces, they have to be put in
the right context. For this, we still refer to the proof of \sp
Theorem \ref{main result} below.                               \end{remarks}

In view of our examples and remarks above, the main importance of
Theorem \ref{main result} below lies in the idea of its proof, here presented
as clearly as possible, free from e.g.\ blurring differential geometric
technicalities. If one wants to get definitely convinced of the theorems
\cite[5.2, 5.3, p.\ 45]{Yam79} not being applicable to a particular (say)
differential geometric problem possibly involving a partial differential
equation, then one should use the proof of Theorem \ref{main result} as a
model, and begin to write a proof of length for example some ten pages.


  \binsubsubhead B{The basic concepts and the main result}\label{subsec B}

Since in \cite{Yam79} various loose notational conventions are utilized making
matters obscure, we first recall the facts from \cite{Yam79} needed below,
reformulated so as to be accordant with the set theoretic notational system we
followed in \cite{Hic} and \cite{SeBGN}\ssp.

\begin{definitions}

For \œ$\sp E\in\LCS(\bosy R\ssp)\,$, let $\sp\SemiNor E$ be the set of all
continuous seminorms on $\sp E\,$. For \œ$\sp N\inc\SemiNor E\,$, we say that
$\sp N$ {\it determines\ssp} $E\sp$ if{}f for every \œ$\sp U\sn\in\ymp E\sp$
there are \œ$\sp\eps\in\Rep$ and a finite \œ$\sp N\aar 0\inc N\sp$ with \œ$\ssp
\bigcap\ssp\{\,p\inve\image[\,0\,,\eps\,]:p\in N\aar 0\ssp\}\inc U\sp$. We say
that $\ssp\Gamma\sp$ is a {\it calibration\ssp} over $\ssp\bosy E\sp$ if{}f we
have \œ$\,\bosy E\in\LCS(\bosy R\ssp)\,^{\dom\bosy E}$ with \œ$\ssp\Gamma \inc
\prodc\{\ssp(\sp\nu\sp,\SemiNor E\ssp):$ $(\sp\nu\sp,E\ssp)\in\bosy E\sp\,\}\sp
$ such that $\ssp\bigcup\ssp\Gamma\ssp\image\sn\{\sp\nu\sp\}$ determines $\sp
E\sp$ whenever $\sp(\sp\nu\sp,E\ssp)\in\bosy E\,$.         \end{definitions}

One easily sees that $N$ determines a given \œ$E\in\LCS(\bosy R\ssp)$ if{}f 
every \œ$\ssp p\in\SemiNor E$ has some finite \œ$N\aar 0\inc N$ and \œ$\smb M
\in\Rep$ with \œ$\ssp p\fvalue x\le\sup\,\{\,\smb M\,r:\exi q\,
(\sp x\ssp,r\sp)\in q\in N\aar 0\ssp\}$ for all \œ$\sp x\in\vecs E\,$. Note
also that if $\Gamma\sp$ is a calibration over $\ssp\bosy E\,$, then $\Gamma\sp
$ and $\ssp\bosy E$ are necessarily small families. It follows that the speech
for example in \cite[Example 1, p.\ 4]{Yam79} of having a calibration over the
class of all normable spaces does not make sense in our set theory. In
\cite[p.\ 3]{Yam79}\ssp, one refers by the term \q{seminorm map} to the
elements of the product set \œ$\ssp
\prodc\{\ssp(\sp\nu\sp,\SemiNor E\ssp):(\sp\nu\sp,E\ssp)\in\bosy E\sp\,\}\,$
which is empty \linebreak (in our set theory) if $\ssp\bosy E\sp$ is a large
family of locally convex spaces.

Arbitrarily fixing any two\ssp-\ssp element set $I\ar 0\,$, for example taking
as $I\ar 0$ the cardinal number \œ$\sp 2\adot = \{\ssp\emptyset\,,1\adot\ssp\}
= \{\ssp\emptyset\,,\{\emptyset\}\sp\}\,$, for the purpose of this note it
would suffice to consider only calibrations over $\sp\bosy E\sp$ with $
  \dom\bosy E=I\ar 0\,$.

\begin{definitions}[\sp Gateaux differentiability\sp]\label{Def Gat} $\null$ \nopagebreak\vskip.6mm

$\rGatDer_EFf = \{\ssp(\sp x\ssp,\sp\ell\,):E\ssp,F\in\LCS(\bosy R\sp)$ and $
                   f\in\vecs F\KP1^{\dom\sn f}$ and \nopagebreak\vskip.2mm $\nKP{19}
     \dom\sn f\in\Nbh(\sp x\ssp,\taurd E\ssp)$ and $\ssp
       \ell\in\Cal L\,(E\ssp,F\ssp)$ and \nopagebreak\vskip.3mm $\nKP{18.2}
     \all{u\in\vecs E\,,\sp V\sn\in\ymp F}\,\exi{\delta\in\Rep}\,\all{t\in\Re}$ \nopagebreak\vskip.2mm $\null\hfill
     0<|\ssp t\ssp|<\delta\imply(\sp t^{\sp\mminus 1}\sp(\sp
     f\sp\fvalue(\sp x+t\,u\sp)\svs E - f\sp\fvalue x\sp)\svs F - \sp
     \ell\ssp\fvalue u\sp)\svs F\in V\sp\,\}\,$, \KP4 \vskip.5mm $\nKP{7.1}
 \tilde f\ssp'(x) = \bigcap\ssp\{\ssp\rGatDer_EFf\sp\fvalue x:
                     \tilde f=(E\ssp,F\sp,f\ssp)\ssp\}\,$. \vskip.5mm

\noin By a real {\it Gateaux differentiable\ssp} map we understand any \œ$\sp
\tilde f=(E\ssp,F\sp,f\ssp)$ such that $E\ssp,F\in\LCS(\bosy R\sp)\sp$ and
$\ssp f\in\vecs F\KP1^{\dom\sn f}\sp$ and $\ssp
          \dom\sn f \inc \dom\sn\rGatDer_EFf\sp$.          \end{definitions}

It follows that if $\sp(E\ssp,F\sp,f\ssp)$ is Gateaux differentiable, then \œ$
\dom\sn f\in\taurd E\sp$ since for every \œ$\sp x\in\dom\sn f\sp$ we have \œ$
 \dom\sn f\in\Nbh(\sp x\ssp,\taurd E\ssp)\,$. For all $\ssp E\ssp,S\sp$
generally hav- ing \œ$\,\sp\roman{Of\KP1}E\subtext{sub\ssp tvs\,}S =
E_{\sp/\ssp S} = (\sp\sigrd E_{\,|\ssp S}\ssp,\taurd E\sp\lei S\ssp)\,$, \,we
consider the Fr\'echet space \vskip.5mm\centerline{$
\Cperinfty(\Re\ssp) = \roman{Of\KP1}\Cinfty(\Re\ssp)\subtext{sub\ssp tvs\,}
  \{\,x:\all{s\in\Re}\,x\fvalue s=x\fvalue(\sp s+1\sp)\ssp\}$} \vskip.5mm

\noin of smooth $\sp 1\ssp$--\,periodic functions
      $\sp\Re\to\Re\ssp$ in the following

\begin{lemma}\label{Le1}

Let \œ$\ssp E=\Cperinfty(\Re\ssp)$ and \œ$\,
f=\seq{\,\varphi\circ x:x\in\vecs E\sp\,}\sp$ where $\ssp\varphi$ is a smooth
function $\ssp\Re\to\Re\sp\,$. Then $\ssp(E\ssp,E\ssp,f\ssp)\ssp'(x) =
\seq{\,\varphi\ssp'\snn\circ x\cdotn v:v\in\vecs E\sp\,}\ssp$ for all $\,
x\in\vecs E\,$.                                                  \end{lemma}

\begin{proof} Let \œ$\ssp
\varphi\ar 1=\{\ssp(\sp s\ssp,t\ssp,\varphi\fvalue t\sp):s\ssp,t\in\Re\sp\,\}\sp
$ and \œ$\sp f\aar 1=\seq{\,\varphi\circ x:x\in\vecs F\sp\,}\,$, where \linebreak
\œ$\sp F=\Cinfty(\Re\ssp)\,$, and first consider the map \œ$\sp\tilde f\aar 1=
(\spp F\sp,F\sp,f\aar 1)\,$. Since we can decompose it \œ$\sp
F\to\Cinfty(\Re\sp\times\Re\ssp)\sqcap F\to F\sp$ by \œ$\,
x\mapsto(\sp\varphi\ar 1\ssp,\spp x\sp)\mapsto\varphi\ar 1\snn\circ[\,\id,x\,]
=\varphi\circ x\ssp$, where the first factor is a continuous affine map, hence
smooth, and the second is smooth by \cite[Theorem 3.6, p.\ 17]{Hic}\ssp, by
the chain rule \cite[Proposition 0.11 p.\ 240]{HiSM} we get \œ$\sp
\tilde f\aar 1\in\CinftyPi(\bosy R\ssp)\,$. Since $\sp E\sp$ is a
(sequentially) closed topological linear subspace of $\ssp F\sp$, and since we
have \œ$\sp f=f\aar 1\ssp|\,\vecs E\sp$ with \œ$\sp\rng f\inc\vecs E\,$, the
assertion of the lemma follows from
\cite[Proposition 3.1, Remarks 3.7\sp(b)\ssp, pp.\ 15, 17\,--\,18]{Hic} in
conjunction with elementary set theoretic manipulations applied to
Definitions \ref{Def Gat} above.                                 \end{proof}

\begin{remarks}

Given a calibration $\ssp\Gamma$ over $\bosy E$ with \œ$(\ssp\mu\,,E\ssp)\ssp,
(\sp\nu\sp,F\ssp)\in\bosy E$ and \œ$L=$ \linebreak $
\nLinb_\Gamma(\ssp\mu\,,\sp\nu\sp)\sbi{\bosy E}\,$, then $L$ is the unique
normable locally convex space with $\sp\sigrd L$ a \nolinebreak vec- \Newline
tor substructure of $\sp\sigrd F\expnota^\sp\upsilon_s\sp E\sp]_{vs}$ and $
\vecs L=\Nu\ssp\inve\image\snn\Repp$ and $\ssp\Nu\ssp\inve\image[\,0\,,1\,]\in
\rajou L\sp\cap\sp\ymp L$ \Newline when we let \œ$\sp
 \Nu = \seq{\,\sup\,\{\,\bosy p\fvalue\nu\fvalue(\ssp\ell\ssp\fvalue v\sp) :
 \bosy p\in\Gamma$ and $\bosy p\fvalue\mu\fvalue v\le 1\ssp\} : \ell \in
 \Cal L\,(E\ssp,F\ssp)\,}\,$. Here we call $\sp\Nu\sp\,|\,\vecs L\sp$ the {\it
canonical $\,\Gamma\,$--\,norm\ssp} \œ$\mu\to\nu$ over $\sp\bosy E\,$. In
\cite{Yam79}\ssp, the space $L$ is imprecisely denoted by \q{\ssp$
\roman{L\sbi{B\Gamma}}(E,F)\sp$}\sp. If also $F$ is sequentially complete,
then $L$ is Banachable, and $\ssp(\sp\sigrd L\,,\sp\Nu\sp\,|\,\vecs L\sp)$ is
a (normed) Banach space, cf.\ \cite[2.3, p.\ 6]{Yam79}\ssp.    \end{remarks}

In view of \cite[5.2, 5.3, p.\ 45]{Yam79}\ssp, the \q{continuous B$\Gamma\,
$--\,differentiability\sp} (or being $\sp\roman{C_{\,B\Gamma}}$) should be a
most important concept in \cite{Yam79}\ssp. Despite this, on page 23 there,
its definition is only vaguely sketched, and we replace this concept by the
generally weaker one given in the following

\begin{definition}\label{defi C_{BGamma}}

A real Gateaux differentiable map \œ$\tilde f=(E\ssp,F\sp,f\ssp)$ we say to be
{\it continuously\ssp} cb $\,\Gamma\,${\it--\,differentiable\ssp} \œ$\mu\to\nu
$ within $\ssp\bosy E\sp$ if{}f $\,\Gamma$ is a calibration over $\sp\bosy E\sp
$ with \œ$(\ssp\mu\,,E\ssp)\ssp,(\sp\nu\sp,F\ssp)\in\bosy E$ such that for \œ$
L=\nLinb_\Gamma(\ssp\mu\,,\sp\nu\sp)\sbi{\bosy E}$ and for $\sp\Nu\sp$ the
canonical $\sp\Gamma\,$-- norm \œ$\ssp\mu\to\nu$ over $\sp\bosy E\,$, we have
\œ$\sp\rng\rGatDer_EFf\inc\vecs L$ and also for \œ$\sp x\in\dom\sn f$ and \linebreak
$\bosy p\in\Gamma$ and $\sp\eps\in\Rep$ there is $\sp\delta\in\Rep$ such that
for all $\sp u\ssp,y\sp$ we have the implication \nopagebreak\vskip.3mm\centerline{$
    \bosy p\fvalue\mu\fvalue u<\delta\,$ and $\,y=x+u\in\dom\sn f\,\imply\,
     \Nu\sp\fvalue(\sp\tilde f\ssp'(\sp y)-\tilde f\ssp'(x))<\eps\,$.}
  \end{definition}

It follows from \cite[2.6, p.\ 29]{Yam79} that if $\sp\dom\sn f\sp$ is convex
or if $\,\Gamma\sp$ is such that for every \œ$\ssp\bosy p\in\Gamma\sp$ and \œ$\ssp
x\in\dom\sn f\sp$ there is \œ$\sp\delta\in\Rep$ with \œ$\sp
\{\,x+u:\bosy p\fvalue\mu\fvalue u<\delta\sp\,\}\inc\dom\sn f\sp$,
cf.\ \cite[pp.\ 18\,--\,19]{Yam79}\ssp, then $\sp\tilde f\sp$ is  continuously
cb $\,\Gamma\,$--\,differentiable \œ$\sp\mu\to\nu\sp$ within $\ssp\bosy E\sp$ \linebreak
if{}f it is \q{\ssp$\roman{C\sbi{B\Gamma}}$}. Consequently, noting that $\sp E\sp
$ is locally convex, if $\sp\tilde f\sp$ is continuously cb $\,\Gamma\,
$--\,differentiable \œ$\sp\mu\to\nu\sp$ within $\ssp\bosy E\,$, for every \œ$\sp
x\in\dom\sn f\sp$ there is $\sp U$ with \œ$\sp x\in U$ such that $\sp
(E\ssp,F\sp,f\,|\,U\ssp)\sp$ is \q{\ssp$\roman{C\sbi{B\Gamma}}$}.
In Remarks \ref{Rems2}\ssp(c) above, we gave the basic ingredi- ents for the
proof of \sp\cite[5.2, p.\ 45]{Yam79} not being applicable to the map $\sp
(\spp G\sp,G\sp,h\,|\,W\ssp)\sp$ anyhow one chooses $\sp W$ with $\sp
   z\ar 0\in W\in\taurd G\ssp$.

\begin{theorem}\label{main result}

Let \œ$\ssp E=\Cperinfty(\Re\ssp)$ and \œ$\,
f=\seq{\,\varphi\circ x:x\in\vecs E\sp\,}$ where $\ssp\varphi$ is a smooth
function \œ$\ssp\Re\to\Re\sp\,$. If
                                    $\ssp(E\ssp,E\ssp,f\ssp)$ is continuously
cb $\,\Gamma\,$--\,differentiable \œ$\ssp\mu\to\nu$ within $\ssp\bosy E\,,$
then there are \œ$\alpha\ssp,\beta\in\Re$ with the property that $\ssp
\varphi=\seq{\,\alpha\,t+\beta:t\in\Re\ssp\,}\,$.              \end{theorem}

\begin{proof} Let \œ$\tilde f=(E\ssp,E\ssp,f\ssp)\,$. Under the premise,
arbitrarily fixing \œ$s\ar 0\in\Re\sp\,$, it suffices to prove indirectly that
\œ$\varphi\ssp''\fvalue s\ar 0=0\,$. To get a contradiction, let \œ$
\varphi\ssp''\fvalue s\ar 0\not=0\,$ and consider \œ$
x=\Re\times\sn\{\ssp s\ar 0\}\,$. Let $\sp\Nu\sp$ be the canonical $\sp
\Gamma\,$-- norm \œ$\ssp\mu\to\nu$ over $\sp\bosy E\,$, and put \œ$
\smb M=\Nu\sp\fvalue(\sp\tilde f\ssp'(x))+1\ssp$. Since $\ssp
\bigcup\ssp\Gamma\sp\image\sn\{\ssp\mu\sp\}$ determines $E\,$, we have \œ$\sp
\Gamma\not=\emptyset\,$, and so we can pick \œ$\sp\bosy p\ar 0\in\Gamma\ssp$.
Taking this $\sp\bosy p\ar 0$ in place of $\sp\bosy p$ and \œ$\sp\eps=1$ in
Definition \ref{defi C_{BGamma}} above, \linebreak there is \œ$\sp
\delta\in\Rep$ such that for \œ$\sp\bosy p\ar 0\KN1\fvalue\mu\fvalue u<\delta$
we have \œ$\sp\Nu\sp\fvalue(\sp\tilde f\ssp'(\sp x+u\sp)-\tilde f\ssp'(x))<1\ssp
$, \linebreak whence further \œ$\sp\Nu\sp\fvalue(\sp\tilde f\ssp'(\sp x+u\sp))
    \le \Nu\sp\fvalue(\sp\tilde f\ssp'(x))+1=\smb M\ssp$. Consequently,
noting also Lemma \ref{Le1} above, for all $\sp\bosy p\in\Gamma\sp$ we have \vskip.4mm \nKP{-2.5}
 (1) \hfill $
 \bosy p\fvalue\nu\fvalue(\sp\varphi\ssp'\snn\circ(\sp x+u\sp)\cdotn v\sp)
  \le\smb M\,(\ssp\bosy p\fvalue\mu\fvalue v\sp)\,$ when $\sp
   \bosy p\ar 0\KN1\fvalue\mu\fvalue u<\delta\,$ and $\,v\in\vecs E\,$. \hfill { \ \ } \vskip.4mm

Take \œ$\sp v=\Re\times\sn\{\sp 1\sp\}\,$. Letting \œ$\,\|\ssp z\ssp\|\sbi i =
  \sup\,\{\,|\,\roman D^{\,l}\sp z\fvalue s\,| :
      l\in i\ssp\yplus\text{ and }\sp s\in\Re\sp\,\}\,$, then also \œ$
\{\ssp\seq{\,\sp\|\ssp z\ssp\|\sbi i:z\in\vecs E\sp\,}:i\in\No\ssp\}$
determines $E\,$. Hence, there is an even \œ$i\ar 0\in\N$ such that we have
the implication \vskip.4mm \nKP{-2.5}
 (2) \hfill $i\ar 0\,\|\ssp u\ssp\|\sbi{i_0}<1\,\imply\,
      \bosy p\ar 0\KN1\fvalue\mu\fvalue u<\delta\,$ for all $\,u\in\vecs E\,$. \hfill { \ \ } \vskip.4mm

\noin Letting \œ$\sp i=i\ar 0\yplus\sp$, we have $i$ odd. Now, there are a
finite \œ$P\inc\Gamma$ and \œ$\smb M\ar 1\in\Rep$ such that for \hfill
$ q = \seq{\,\sup\,\{\,\smb M\ar 1\ssp r:\exi{\ssp\bosy p\in P}\,
       (\sp z\ssp,r\sp)\in\bosy p\fvalue\nu\,\}:z\in\vecs E\sp\,}$ \hfill
and \vskip.3mm $\nKP{22.3}
  p = \seq{\,\sup\,\{\,\smb M\,\smb M\ar 1\ssp r:\exi{\ssp\bosy p\in P}\,
       (\sp z\ssp,r\sp)\in\bosy p\fvalue\mu\,\}:z\in\vecs E\sp\,}\,$, \vskip.3mm

\noin we have \œ$\sp\|\ssp z\ssp\|\sbi i\le q\fvalue z$ for all \œ$z\in\vecs E\,
$. Having \œ$p\in\SemiNor E\,$, there further is \œ$\sp j\in\No$ \linebreak
with \œ$\ssp p\fvalue z\le j\,\|\ssp z\ssp\|\sbi j$ for all \œ$z\in\vecs E\,$.
Using these and \œ$\sp j\,\|\ssp v\ssp\|\sbi j=j\,1=j\sp\ydot\sp$, from (1)
and (2) we obtain \vskip.4mm \nKP{-2.5}
 (3) \hfill $\|\,\varphi\ssp'\snn\circ(\sp x+u\sp)\,\|\sbi i\le j\sp\ydot\sp$
 for all $\,u\in\vecs E\sp$ with $\ssp i\ar 0\,\|\ssp u\ssp\|\sbi{i_0}<1\ssp$. \hfill { \ \ } \vskip.4mm

Next, noting that \œ$\ssp|\,r+t\,|\ge|\ssp r\ssp|-|\ssp t\ssp|\sp$ for \œ$\sp
r\sp,t\in\Re\sp\,$, and utilizing the quite combinatorial idea in the proof of
\cite[Proposition 10, pp.\ 6\,--\,7]{Colo}\ssp, one deduces existence of $\sp
\smb M\ar 2\in\Rep$ such that for all $\sp u\in\vecs E\sp$ and $\sp s\in\Re\sp$
with $\ssp\|\ssp u\ssp\|\sbi{0.}\le 1\sp$ we have \vskip.4mm \nKP{-2.5}
 (4) \hfill $
    |\,\roman D^{\ssp i}\sp(\sp\varphi\ssp'\snn\circ(\sp x+u\sp))\fvalue s\,|
\ge |\,\varphi\ssp''\snn\circ(\sp x+u\sp)\snn\cdot
    \roman D^{\ssp i}\sp u\fvalue s\,| - \smb M\ar 2\,
     (\ssp 1 + \|\ssp u\ssp\|\sbi{i_0})^{\ssp i}\,$. \hfill { \ \ } \vskip.4mm

\noin With \œ$\smb A=|\,\varphi\ssp''\fvalue s\ar 0\,|\in\Rep\sp$, now
choosing \œ$n\in\Zepp$ so that we have the inequalities $ j\sp\ydot <
\smb A\,(\sp 2\,\pi\,n\sp)^{\,\frac 12} - \smb M\ar 2\,\big(\ssp 1 +
(\sp 2\,\pi\,n\sp)^{\sp\mminus\frac 12}\sp\big)\RHB{.2}{\sp^i}\sp$ and $\,
i\ar 0\ssp(\sp 2\,\pi\,n\sp)^{\sp\mminus\frac 12} < 1\ssp$, we take \vskip.3mm\centerline{$
 u = \Seq{\,(\sp 2\,\pi\,n\sp)^{\sp\mminus i_0\sn\yydot - \ssp\frac 12}\sp
      \sin\ssp(\sp 2\,\pi\,n\,(\sp s-s\ar 0)):s\in\Re\sp\,}\,$.} \vskip.3mm

\noin Then we have \œ$\ssp i\ar 0\,\|\ssp u\ssp\|\sbi{i_0} =
i\ar 0\ssp(\sp 2\,\pi\,n\sp)^{\sp\mminus\frac 12} < 1\ssp$, whence recalling
that $i\sp$ is odd, by (3) and (4) we obtain \vskip.5mm $\nKP4
 j\sp\ydot < \smb A\,(\sp 2\,\pi\,n\sp)^{\,\frac 12} - \smb M\ar 2\,
 \big(\ssp 1 + (\sp 2\,\pi\,n\sp)^{\sp\mminus\frac 12}\sp\big)\RHB{.2}{\sp^i}$ \nopagebreak\vskip.3mm $\nKP8
 = |\,\varphi\ssp''\fvalue s\ar 0\,|\,
   |\,\roman D^{\ssp i}\sp u\fvalue s\ar 0\,| - \smb M\ar 2\,
              (\ssp 1 + \|\ssp u\ssp\|\sbi{i_0})^{\ssp i}$ \nopagebreak\vskip.3mm $\nKP8
 = |\,\varphi\ssp''\snn\circ(\sp x+u\sp)\snn\cdot
       \roman D^{\ssp i}\sp u\fvalue s\ar 0\,| - \smb M\ar 2\,
              (\ssp 1 + \|\ssp u\ssp\|\sbi{i_0})^{\ssp i}$ \nopagebreak\vskip.3mm $\nKP8
\le |\,\roman D^{\ssp i}\sp(\sp\varphi\ssp'\snn\circ
                       (\sp x+u\sp))\fvalue s\ar 0\,|
\le \|\,\varphi\ssp'\snn\circ(\sp x+u\sp)\,\|\sbi i
\le j\sp\ydot\sp$, \,a contradiction.                            \end{proof}

Note that basic idea in the preceding proof is the same which we already
uti- lized when establishing \cite[Propositions 3\ssp, 5]{Muller}\ssp.

\begin{remark}

Fixing a calibration is required also for the inFT \cite[8.4, p.\ 457]{Y-Bull}
but this result has a nature quite different from that of the IFTs in
\cite[p.\ 45]{Yam79}\ssp. It seems that \cite[8.4]{Y-Bull} has genuine
applications although quite an amount of work is required for the verification
of the premise. Here, by a {\it genuine\ssp} application of an IFT
generalizing the corresponding Banach space theorem we mean an application not
covered by the classical theorem.                               \end{remark}


\end{document}